\documentclass[12pt]{article}
\usepackage{amsmath}
\usepackage{amsmath}
\usepackage{amssymb}
\usepackage{dsfont}
\usepackage{cite}
\newsymbol \blackbox 1004
\newcommand{\eh}{\hfill}\newlength{\sperr}

\newenvironment{proof}{{\settowidth{\sperr}{\bf\rm
Proof}%
\par\addvspace{0.3cm}\noindent\parbox[t]{1.3\sperr}
{\bf\rm P\eh r\eh o\eh o\eh f\eh }%
}}{\nopagebreak\mbox{}
$\blackbox$\par\addvspace{0.3cm}}

\def\nn{\nonumber}

\def\a{\alpha}
\def\b{\beta}
\def\bt{\beta}
\def\g{\gamma}

\def\vk{\varkappa}

\def\Lam{\Lambda}

\def\la{\lambda}
\def\om{\omega}

\def\Th{\Theta}
\def\ze{\zeta}
\def\Up{\Upsilon}
\def\vp{\varphi}
\def\vt{\vartheta}
\def\ve{\varepsilon}
\def\wh{\widehat}
\def\wt{\widetilde}
\def\ov{\overline}

\def\p{\partial}

\def\BC{{\mathbb C}}

\def\BR{{\mathbb R}}
\def\BN{{\mathbb N}}
\def\clp{{\mathcal P}}
\def\cla{{\mathcal A}}
\def\clb{{\mathcal B}}

\def\cle{{\mathcal E}}
\def\clf{{\mathcal F}}

\def\clh{{\mathcal H}}
\def\clj{{\mathcal J}}

\def\cli{{\mathcal I}}
\def\clk{{\mathcal K}}

\def\cln{{\mathcal N}}

\def\clw{{\mathcal W}}

\def\clr{\mathcal{R}}
\def\cls{\mathcal{S}}

\def\clv{{\mathcal V}}
\def\cly{{\mathcal Y}}
\def\clz{{\mathcal Z}}
\def\cld{{\mathcal D}}

\def\diag{\mathrm{diag}}
\newcommand{\E}{\mathrm{e}}
\newcommand{\I}{\mathrm{i}}
\def\mf{\mathfrak}

\def\bT{{\bf T}}
\def\bB{{\bf B}}

\newtheorem{Pa}{Paper}[section]
\newtheorem{Tm}[Pa]{{\bf Theorem}}

\newtheorem{Cy}[Pa]{{\bf Corollary}}
\newtheorem{Rk}[Pa]{{\bf Remark}}

\newtheorem{Ee}[Pa]{{\bf Example}}
\newtheorem{Dn}[Pa]{{\bf Definition}}

\newtheorem{Pn}[Pa]{{\bf Proposition}}

\title{On a class of  canonical systems corresponding to matrix string equations: general-type and explicit fundamental
solutions  and  Weyl--Titchmarsh theory}

\author{Alexander Sakhnovich}

\date{}
\begin{document}
\maketitle

\begin{abstract}  An important representation of the general-type fundamental solutions
of the canonical systems corresponding to matrix string equations is established
using linear similarity of a certain class of Volterra operators to the squared integration.
Explicit fundamental solutions of these canonical systems are also constructed
via the GBDT version of Darboux transformation. Examples and applications to
dynamical canonical systems are given. Explicit solutions of the dynamical canonical systems are constructed as well.
Three appendices are
dedicated to the Weyl--Titchmarsh theory for canonical systems,
transformation of a subclass of canonical systems into matrix  string equations (and of a smaller subclass of canonical systems into matrix Schr\"odinger equations), and a linear similarity problem for Volterra operators. 
\end{abstract}


\vspace{0.2em}


\section{Introduction} \label{intro}
\setcounter{equation}{0}
Canonical (spectral canonical) systems have the form
\begin{align} &       \label{1.1}
w^{\prime}(x,\la)=\I \la J H(x)w(x,\la), \quad J:=\begin{bmatrix} 0 & I_p \\ I_p & 0\end{bmatrix} \quad \Big(w^{\prime}:=\frac{d}{dx}w\Big),
 \end{align} 
where $\I$ is the  imaginary unit ($\I^2=-1$), $\la$ is the so called spectral parameter,
$I_{p}$ is the $p \times p$ $(p\in \BN)$ identity
matrix,  $\BN$ stands for the set of positive integer numbers,  $H(x)$ is a $2p \times 2p$ matrix valued function (matrix function),
and 
$H(x) \geq 0$ (that is, the matrices $H(x)$ are self-adjoint and the eigenvalues of $H(x)$ are nonnegative).
Canonical systems are important objects of analysis, being perhaps the most important class of  the one-dimensional
Hamiltonian systems and including (as subclasses) several classical equations. They have been actively studied in many
already classical as well as in various recent works (see, e.g., \cite{ArD, dBr,  EKT, FKS, GoKr, Mog, Rem, Rom, RW, Rov, ALS2019, SaSaR, SaL2, Su} 
and numerous references
therein). We will also consider (and construct explicit solutions) for a more general class
 \begin{align} &       \label{1.2}
 w^{\prime}(x,\la)=\I \la j H(x)w(x,\la), \quad H(x)=H(x)^*, \quad j:=\begin{bmatrix} I_{m_1} & 0 \\  0 & -I_{m_2}\end{bmatrix},
 \end{align} 
where  $\, m_1,m_2 \in \BN$. Here, we set 
$$m_1+m_2=:m,$$
$H$ is an $m\times m$ locally integrable matrix function, and $H(x)^*$ means the complex conjugate transpose of the matrix $H(x)$.
System \eqref{1.2} will be called a {\it generalized canonical system} and the corresponding matrix function $H$ will be called 
a {\it generalized Hamiltonian}.

In the case $m_1=m_2=:p$, it is easily checked that $j$ and $J$ are unitarily similar:
\begin{align}& \label{1.3-}
J=\Theta j \Theta^*, \quad \Theta:=\frac{1}{\sqrt{2}}\begin{bmatrix} I_p & -I_p \\ I_p & I_p\end{bmatrix} ,
\end{align}
that is (assuming $H(x)\geq 0$), system \eqref{1.2} is equivalent to \eqref{1.1} (see Appendix \ref{String} for details).  We call the system
\begin{align}& \label{1.3}
 w^{\prime}(x,\la)=\I \la j H(x)w(x,\la), \quad H(x)\geq 0 \quad (m_1=m_2=p)
\end{align}
canonical (as well as the equivalent system \eqref{1.1}). The matrix function $H(x)$ is called the Hamiltonian of this system.

In most works on canonical systems the  less complicated $2 \times 2$ Hamiltonian case (i.e, the case $p=1$) is dealt with
although the cases with other values of $p$ ($p>1$) are equally important. Interesting recent works \cite{EK, Langer, Wor}  on $2\times 2$
canonical systems and string equations  also contain some useful references.
Here, we deal with the case of $2p \times 2p$ Hamiltonians $(p\geq 1)$.

Well-known Dirac (or Dirac-type) systems are equivalent
to a special subclass of canonical systems  (see \cite{FKS, GeS, ALS2019, SaSaR, SaL2} and references therein). The Hamiltonians corresponding to Dirac
systems (after we switch from the representation \eqref{1.1} to the representation \eqref{1.3}),
have the form
\begin{align}& \label{C5}
 H(x)=\g(x)^*\g(x), \quad \g(x)j \g(x)^*=-I_p,
\end{align}
where $\g$ are $p\times 2p$ matrix functions.
For instance,  formulas (1.7), (1.11), (1.12), and (1.26) in \cite{SaSaR}  lead to the representation \eqref{C5}.

The Hamiltonians, which we consider in this paper, have the form
\begin{align}& \label{I1}
 H(x)=\b(x)^*\b(x), \quad \b(x)j \b(x)^*=0 \quad (p \geq 1),
\end{align}
where $\b$ are again  $p\times 2p$ matrix functions. Thus, canonical systems \eqref{1.3} with Hamiltonians of the form
\eqref{I1} are dual in a certain way to the class of canonical systems corresponding to Dirac systems.
Under some natural conditions, systems \eqref{1.3}, \eqref{I1} are also equivalent to the matrix string equations (see \cite[Chapter 11]{SaL2} and Appendix \ref{String} in our paper).
Canonical systems \eqref{1.3} with some special Hamiltonians of the form \eqref{I1} appear, for instance, as the  linear systems auxiliary to {\it nonlinear
second harmonic generation equations} \cite{Kau, KauS}.

We note that first order symplectic systems $y^{\prime}(x)=F(x)y(x)$ are actively studied (see \cite{BoDo, DoO, DoH, KraS} and references therein).
In the reformulation for our case, simplecticity means the equality
$$F(x)^*j+jF(x)+\mu(x)F(x)^*jF(x)=0.$$
Thus, generalized canonical systems \eqref{1.2} (where $m_1=m_2$) are symplectic, with $\mu \equiv 0$.
Canonical systems \eqref{1.3}, \eqref{I1}, which are our main topic in this paper, remain symplectic for any choice
of $\mu(x)$.

The normalization condition
\begin{align}& \label{I1+}
\b^{\prime}(x)j\b(x)^*=\I I_p
\end{align}
for Hamiltonians of the form \eqref{I1} is essential for the construction of fundamental solutions
and solving inverse problems. In  Appendix \ref{String}, we show that matrix Schr\"odinger equations
may be transformed into canonical systems \eqref{1.3}, \eqref{I1}, \eqref{I1+} satisfying certain additional
condition. There is considerable interest in generalized Schr\"odinger equations (e.g., in Schr\"odinger equations
with distributional potentials, see some references in \cite{EGNST}).
One can say that  systems \eqref{1.3}, \eqref{I1}, \eqref{I1+} present an important generalization of the
of the matrix Schr\"odinger equations.
Canonical systems with Hamiltonians satisfying \eqref{I1}, \eqref{I1+}
were briefly considered in \cite{LA94, SaL2}. However, local boundedness of $\b^{\prime\prime}$ was required
there instead of the local  square-integrability of  $\b^{\prime\prime}$, which we require in the next section.
In Section~\ref{FuSo}, we represent the
fundamental solutions for this case as the transfer matrix function  from \cite{SaL1, LA94, SaL2}.
For this purpose, we use the linear similarity
of the operator $K=\I \b(x)j\int_0^x\b(t)^*\cdot dt$ to the operator \eqref{AC2}  of  squared integration
as well as the form of the corresponding similarity transformation operator $V$ (see Theorem \ref{TmSim} and its proof in Appendix \ref{Sim}).

The representation of the fundamental solutions in Section \ref{FuSo} is important in itself and (in view of the interconnections
between fundamental solutions and Weyl--Titchmarsh functions) it also presents 
a crucial step in solving the inverse problem to recover canonical system from the Weyl--Titchmarsh function.

Some basic results and notions on the Weyl-Titchmarsh theory of the general-type canonical systems
\eqref{1.3} are described in Appendix \ref{Weyl}.
The results are conveniently reformulated  in terms of system \eqref{1.3} instead
of system \eqref{1.1}, and, what is essentially more important, certain redundant
conditions contained in \cite[Appendix A]{SaSaR} are removed.

In other sections of the paper we study explicit solutions of systems \eqref{1.2}
with generalized Hamiltonians $dj+\b(x)^*\b(x)$ as well as explicit solutions
and corresponding Weyl--Titchmarsh (Weyl) functions of the canonical systems \eqref{1.3}, \eqref{I1}.
We note that explicit solutions of Dirac systems and the corresponding Weyl--Titchmarsh theory
have been studied sufficiently well (see, e.g., \cite{GKS2, SaA94, SaSaR})
but the situation with the systems \eqref{1.3}, \eqref{I1} is quite different.

Explicit solutions of canonical systems and their properties are of essential theoretical and applied interest.
Various versions of B\"acklund-Darboux transformations and related dressing and 
commutation methods \cite{Ci, CoIv, Ge, GeT,  Gu, KoSaTe,
Mar, MS, ZM} are fruitful tools in the construction of explicit solutions of linear and integrable nonlinear equations.
B\"acklund-Darboux transformations for canonical and dynamical canonical systems, respectively, were constructed
in \cite{SaA97} and \cite{ALS17}. More precisely,  GBDT (generalized B\"acklund-Darboux transformation)
was constructed for these systems. It is important that GBDT (see, e.g., \cite{GKS2, KoSaTe, SaA94, SaSaR, ALSgrav} and references
therein) is characterized by the generalized matrix eigenvalues (not necessarily diagonal) and
the corresponding generalized  eigenfunctions. In  Section \ref{GBDT},  generalized matrix eigenvalues  and
 generalized  eigenfunctions are denoted by $\cla$ and $\Lam(x)$, respectively. 
 
 Although GBDT 
 for canonical systems was obtained in \cite{SaA97}, a crucial  step of constructing  the generalized eigenfunctions
 $\Lam(x)$ (which is necessary for constructing explicitly Hamiltonians and fundamental solutions)
 is done in the present paper. 
 More precisely, {\it the procedure works in the following way}. We start with some {\it initial systems}  \eqref{1.2}, where
 {\it initial Hamiltonians} $H(x)$ are comparatively simple,
 and construct explicitly
 the fundamental
 solutions and generalized eigenfunctions for these systems. (In particular, some considerations from \cite{SaA05, ALSgrav}
 were helpful for this purpose.)
 Using generalized eigenfunctions, the {\it transformed generalized Hamiltonians} and so called {\it Darboux matrices} are constructed as well. 
 Recall that Darboux matrix for generalized
 canonical systems is the matrix 
 function $\Up(x,\la)$ satisfying the equation 
 $$\Up^{\prime}(x,\la)=\I\la\big(j\wt H(x)\Up(x,\la)- \Up(x,\la)jH(x)\big),$$
 where $H$ is the initial generalized Hamiltonian and $\wt H$ the transformed one. In this way, we obtain
 fundamental solutions $\wt W$ for a wide class of the {\it transformed systems} (i.e., systems with the transformed
 generalized Hamiltonians $\wt H(x))$. Indeed, it is easy to see that $\wt W$ is expressed via the
 fundamental solution $W$ of the initial system and the Darboux matrix, namely, $\wt W(x,\la)=\Up(x,\la)W(x,\la)$.
 
 Some preliminaries on GBDT for the generalized canonical systems are given, and transformed generalized Hamiltonians and Darboux matrices
 are constructed in Section \ref{GBDT}.  Generalized eigenfunctions are constructed explicitly in Section \ref{Tran}.
Explicit formulas for fundamental solutions of the initial systems  and for the Weyl functions of the transformed canonical systems on the semi-axis $[0,\infty)$
 are established in Section \ref{Can}. It is shown in Section~\ref{MSchr} that  the second equality in \eqref{I1} (i.e., the equality $\b(x) j\b(x)^*=0$) for the
 initial matrix function $\b(x)$
 yields the equalities $\wt \b(x) j \wt \b(x)^*=0$ and $\wt \b(x)^{\prime}j \wt\b(x)^*=\b(x)^{\prime}j\b(x)^*$ for the transformed
 matrix function $\wt \b(x)$.
Some interesting examples
 are treated in Section \ref{ExAp}.
 
 There are important connections between spectral and dynamical characteristics as well as between spectral and dynamical systems
(see, e.g., \cite{Beli, BreS, JLP, ALS15, ALS17, SaL18} and references therein).  In particular, GBDT for the spectral canonical
systems  \eqref{1.3} is closely related to the GBDT for the dynamical canonical system
\begin{align}& \label{I2}
H(x) \frac{\p}{\p t}Y(x,t)=j\frac{\p}{\p x}Y(x,t) \quad (m_1=m_2=p), \quad H(x)\geq 0, \quad  x\geq 0.
\end{align}
We note that the invertibility of $H(x)$ was assumed for the dynamical canonical system considered in \cite{ALS17}, and 
system \eqref{I2} slightly differs from the one in \cite{ALS17}.
Dynamical canonical systems are of interest in mechanics and control theory (see, e.g., \cite{JZ}).
The GBDT formula for $Y$ and some explicit examples of $H$ and $Y$ are also discussed in Section \ref{ExAp}.

As usual, $\BR$ stands for the real axis, $\BR_+=\{r: \, r\in \BR, \,\, r\geq 0\}$, $\BC$ stands for the complex plane,
the open upper half-plane is denoted
by $\BC_+$, and  $\ov{a}$ means the complex conjugate of $a$. The notation $\Re(a)$ stands for the real part of $a$,
and $\Im(a)$ denotes the imaginary part of $a$.
The notation $\diag\{d_1, \ldots\}$ stands for the
diagonal (or block diagonal) matrix with the entries (or blocks) $d_1, \ldots$ on the main diagonal.
The space of square-integrable functions on $(0,\,b)$ $(0<b\leq \infty)$ is denoted by $L_2(0, b)$ and the corresponding space
of $p$-dimensional column vector functions is denoted by $L_2^p(0, b)$.
By $L_2^{p\times q}(0, b)$ we denote the class of $p\times q$ matrix functions with the entries belonging
to $L_2(0, b)$. The notation $I$ stands for the identity operator. The norm $\|A\|$ of the $n\times n $ matrix $A$
means the norm of $A$ acting in the space $\ell_2^n$ of the sequences of length $n$. 
The class of bounded operators acting from the Hilbert space $\clh_1$ into Hilbert space $\clh_2$ is denoted by
$\bB(\clh_1,\clh_2)$, and we set $\bB(\clh):=\bB(\clh,\clh)$.

\section{General-type fundamental solutions} \label{FuSo}
\setcounter{equation}{0}
In this section, we study canonical system \eqref{1.3} satisfying conditions \eqref{I1} and \eqref{I1+}:
\begin{align} & \label{F1}
H(x)=\b(x)^*\b(x), \quad \b(x)j \b(x)^*=0, \quad \b^{\prime}(x)j\b(x)^*=\I I_p.
\end{align}
Let us consider the system \eqref{1.3}, \eqref{F1} on some finite interval $[0,\bT]$ ($\bT>0$).
The linear similarity of the operators $K\in \bB\big(L_2^p(0,\bT)\big)$ and $A\in \bB\big(L_2^p(0,\bT)\big)$, where
\begin{align} & \label{AC2}
Kf=\I \b(x)j\int_0^x\b(t)^* f(t)dt, \quad Af=\int_0^x(t-x)f(t) dt,
\end{align}
is essential for us. Here,
the operator $A$ is introduced  as  the squared integration
multiplied by $-1$. (Recall that in the case of Dirac systems the analog of $A$ is the integration multiplied by $\I$.)
It is easy to see that
\begin{align} & \label{F2}
K-K^*=\I \b(x)j\int_0^{\bT}\b(t)^* \cdot dt.
\end{align}
If $\b^{\prime\prime}(x)\in L_2^{p\times 2p}(0,\bT)$, we have (according to Theorem \ref{TmSim}) $K=VAV^{-1}$, which we substitute
into \eqref{F2}. Multiplying both parts of the derived equality by $V^{-1}$ from the left  and by  $(V^*)^{-1}$ from the right, 
we obtain the operator identity
\begin{align} & \label{F3}
AS-SA^*=\I\Pi j \Pi^*,
\end{align}
where
\begin{align} & \label{F4}
S=V^{-1}(V^*)^{-1}>0, \quad \Pi h=\Pi(x)h, \quad \Pi(x):=\big(V^{-1}\b\big)(x), \\
& \label{F5}
\Pi\in \bB\big(\BC^{2p},\, L_2^p(0,\bT)\big), \quad  \Pi(x) \in L_2^{p\times 2p}(0,\bT),
\quad h\in \BC^{2p}.
\end{align}
Note that $\Pi$ above is the operator of multiplication by the matrix function $\Pi(x)$ and the operator $V^{-1}$ is applied  to $\b$ (in the expression
$V^{-1}\b$) columnwise. The transfer matrix function corresponding  to the so called $S$-node (i.e., to the triple $\{A,S,\Pi\}$ satisfying \eqref{F3})
has the form
\begin{align} & \label{F6}
w_A(\la)=w_A(\bT,\la)=I_{2p}-\I j\Pi^*S^{-1}(A-\la I)^{-1}\Pi,
\end{align}
and was first introduced and studied in \cite{SaL1}. We introduce the projectors $P_{\ell}\in \bB\big(L_2^p(0,\bT), \, L_2^p(0,\ell)\big)$:
\begin{align} & \label{F7}
\big(P_{\ell}f\big)(x)=f(x) \quad (0 < x < \ell, \quad \ell \leq \bT).
\end{align}
Now, we set
\begin{align} & \label{F8}
S_{\ell}=P_{\ell}SP_{\ell}^*, \quad V_{\ell}=P_{\ell}VP_{\ell}^*, \quad A_{\ell}=P_{\ell}AP_{\ell}^*, \quad \Pi_{\ell}=P_{\ell}\Pi,
\\ & \label{F9}
w_A(\ell,\la)=I_{2p}-\I j\Pi_{\ell}^*S_{\ell}^{-1}(A_{\ell}-\la I)^{-1}\Pi_{\ell}.
\end{align}
Since $V$ is a triangular operator, $V^{-1}$ is triangular as wel, and we have $P_{\ell}V^{-1}= P_{\ell}V^{-1}P_{\ell}^*P_{\ell}$.
Hence, taking into account \eqref{F4} and \eqref{F8} we derive
\begin{align} & \label{F10}
P_{\ell}V^{-1}P_{\ell}^*V_{\ell}=P_{\ell}V^{-1}P_{\ell}^*P_{\ell}VP_{\ell}^*=P_{\ell}V^{-1}VP_{\ell}^*=I,
\\ & \label{F11}
S_{\ell}=P_{\ell}V^{-1}(V^*)^{-1}P_{\ell}^*=P_{\ell}V^{-1}P_{\ell}^*P_{\ell}(V^*)^{-1}P_{\ell}^*.
\end{align}
It follows that
\begin{align} & \label{F12}
V_{\ell}^{-1}=P_{\ell}V^{-1}P_{\ell}^*, \quad S_{\ell}=V_{\ell}^{-1}(V_{\ell}^*)^{-1}.
\end{align}
We also have $P_{\ell}A= P_{\ell}AP_{\ell}^*P_{\ell}$. Thus, multiplying both parts of  \eqref{F3} by $P_{\ell}$ from the left and by  $P_{\ell}^*$ from
the right (and using \eqref{F8}, \eqref{F12}, and the last equality in \eqref{F4}) we obtain
\begin{align} & \label{F13}
A_{\ell}S_{\ell}-S_{\ell}A_{\ell}^*=\I\Pi_{\ell} j \Pi_{\ell}^*, \quad \Pi_{\ell}(x)=\big(V_{\ell}^{-1}\b\big)(x) \quad (0<x<\ell).
\end{align}
Clearly $w_A(\ell,\la)$ coincides with $w_A(\bT, \la)$ when $\ell=\bT$.
\begin{Rk}\label{wA}
Relations \eqref{F9}, \eqref{F12} and \eqref{F13} show that $S_{\ell}$ and $w_A(\ell,\la)$ may be defined via $V_{\ell}$ $($and $\b(x)$ given on $[0,\ell])$
precisely in the same way as 
$w_A(\bT,\la)$ is constructed via $V$ $($and $\b(x)$ given on $[0,\bT])$. Moreover, according to Remark \ref{RkLT}, $V_{\ell}$ may be constructed in the same way as $V$, and so $w_A(\ell,\la)$ does not depend on the choice of $\b(x)$ for $\ell<x<\bT$ and the choice of $\bT\geq \ell$. In particular, $w_A(\ell,\la)$ is uniquely defined on the semi-axis
$0<\ell<\infty$ for $\b(x)$  considered on the semi-axis $0\leq x<\infty$.
\end{Rk}
The fundamental solution of the canonical system \eqref{1.3},
where Hamiltonian has the form \eqref{F1} may be expressed via the  transfer functions $w_A(\ell,\la)$ using continuous
factorization theorem \cite[p. 40]{SaL2} (see also \cite[Theorem 1.20]{SaSaR} as a more convenient
for our purposes presentation).
\begin{Tm} Let the Hamiltonian of the canonical system \eqref{1.3} have the form \eqref{F1}, where $\b(x)$ is a $p\times 2p$ matrix function
two times differentiable and such that $\b^{\prime\prime}(x)\in L_2^{p\times 2p}(0,\bT)$, if the canonical system is considered on 
the finite interval $[0,\bt]$, or the entries of $\b^{\prime\prime}(x)$ are locally square integrable, if the canonical system is considered on $[0,\infty)$.

Then, the fundamental solution $W(x,\la)$ of the canonical system normalized by $W(0,\la)=I_{2p}$ admits representation
\begin{align} & \label{F14}
W(\ell,\la)=w_A\Big(\ell,\frac{1}{\la}\Big).
\end{align}
\end{Tm}
\begin{proof}. First, we fix some $0<\bT<\infty$ and consider $\b(x)$ on $[0,\bT]$. It is easy to
see that the projectors $P_{\ell}$ and the triple $\{A,\, S,\, \Pi\}$ satisfy conditions of \cite[Theorem 1.20]{SaSaR}.
Hence, according \cite[Theorem 1.20]{SaSaR} the matrix function $w_A\Big(\ell,\frac{1}{\la}\Big)$ is the normalized
fundamental solution of the canonical system \eqref{1.3} with Hamiltonian
\begin{align} & \label{F15}
H(\ell)=\frac{d}{d \ell}\int_0^{\ell}\Pi_{\ell}(x)^*S_{\ell}^{-1}\Pi_{\ell}(x)dx,
\end{align}
where $S_{\ell}^{-1}$ is applied to $\Pi_{\ell}(x)$ columnwise. Using the second equalities in \eqref{F12} and \eqref{F13},
we rewrite \eqref{F15} in the form
\begin{align} & \label{F16}
H(\ell)=\frac{d}{d \ell}\int_0^{\ell}\b(x)^*\b(x)dx=\b(\ell)^*\b(\ell),
\end{align}
and the statement of the theorem is proved on $[0,\bT]$. Taking into account Remark \ref{wA}, we see
that the statement of the theorem is valid on $[0,\infty)$ as well.
\end{proof}
\begin{Rk} The operators $\cls_{\ell}$ satisfying \eqref{F13} are so called structured operators.
The study of the structured operators in inverse problems takes roots in the seminal note \cite{Krein} by M.G.  Krein
and was developed by L.A. Sakhnovich in \cite{SaL2-, LA94, SaL2}.
\end{Rk}
\section{GBDT: Darboux matrices for \\ generalized canonical systems} \label{GBDT}
\setcounter{equation}{0}
Let us consider systems \eqref{1.2} on finite or semi-infinite intervals. Without loss of generality, we choose either the
intervals $\cli_{\bT}=[0,\bT]$ $(\bT<\infty)$ or the semi-axis $\BR_+=[0,\infty)$.
We also fix an initial generalized Hamiltonian  $H(x)=H(x)^*$. 
Given an initial $m \times m$ generalized Hamiltonian $H(x)$, each GBDT is (as usual) determined by some $n \times n$ matrices $\cla$ and $\cls(0)=\cls(0)^*$ ($n \in \BN$)
and by an $n \times m$ matrix $\Lam(0)$ which satisfy the matrix identity 
\begin{align}& \label{D8}
\cla \cls(0)-\cls(0)\cla^*=\I \Lam(0) j \Lam(0)^*.
\end{align}
Taking into account the initial values $\Lam(0)$ and $\cls(0)$ (and using the matrix $\cla$ and the matrix function $H(x)$) we introduce
matrix functions $\Lam(x)$ and $\cls(x)=\cls(x)^*$ via the equations:
\begin{align}& \label{D9}
\Lam^{\prime}(x)=-\I \cla \Lam(x)jH(x) , \quad \cls^{\prime}(x)=\Lam(x)jH(x)j\Lam(x)^*.
\end{align}
It is easy to see that \eqref{D8} and \eqref{D9} yield  \cite{SaA97} the identity
\begin{align}& \label{D10}
\cla \cls(x)-\cls(x)\cla^*\equiv \I \Lam(x) j \Lam(x)^*.
\end{align}
\begin{Rk}\label{RkDM}
We note that, similar  to the case of the general-type fundamental solutions in Section \ref{FuSo}, 
we also use operator identities and transfer matrix function in 
Lev Sakhnovich  form in our GBDT constructions. However, instead of the infinite-dimensional operators in Section \ref{FuSo}, identities \eqref{D8} and \eqref{D10}
are written for matrices. Here, we use  calligraphic letter $\cla$ and $\cls$ instead $A$ and $S$ in Section \ref{FuSo} $($and the
notation $\Lam$ instead of $\Pi)$ for the elements of the $S$-node $($of the triple $\{\cla, \cls,\Lam\})$.
The so called Darboux matrix from Darboux transformations  is represented in GBDT $($for each $\, x\, )$  as the transfer matrix function.   
More precisely, we will show that in the points of invertibility of $\cls(x)$ $($for the case of the  generalized canonical system$)$ the Darboux matrix
is expressed via
\begin{align}& \label{D13}
w_\cla(x,\la)=I_m-\I j \Lam(x)^*\cls(x)^{-1}(\cla-\la I_n)^{-1}\Lam(x)
\end{align}
$($see \eqref{D17}$)$. The dependence of $\cls, \Lam$ and $w_{\cla}$ on $x$ is of basic importance and greatly differs
from the dependence of $S_{\ell}$, $\Pi_{\ell}$ and $w_A(\ell,\la)$ on $\ell$ in Section \ref{FuSo}.
\end{Rk}
According to \cite{SaA97}, $w_\cla(x,\la)$ satisfies the equation
\begin{align}& \label{D11}
w_\cla^{\prime}(x,\la)=\big(\I \la j H(x)-\wt q_0(x)\big)w_\cla(x,\la)-\I \la w_\cla(x,\la) jH(x); 
\\ & \label{D12}
 \wt q_0(x):=j\Lam(x)^*\cls(x)^{-1}\Lam(x) j H(x)-j H(x) j \Lam(x)^*\cls(x)^{-1}\Lam(x).
\end{align}
Note that \eqref{D11}  follows directly from 
\eqref{D9}--\eqref{D13}. Moreover, \eqref{D10} yields (see \cite{SaA97} or \cite[(1.88)]{SaSaR}):
\begin{align} \nn
w_\cla(x, \ov{\mu})^*jw_\cla(x,\la)=&j+\I(\mu-\la) 
\\ & \label{D14}
\times \Lam(x)^*(\cla^*-\mu I_n)^{-1}\cls(x)^{-1} (\cla-\la I_n)^{-1}\Lam(x).
\end{align}
Relation \eqref{D15} shows that, under conditions $\det((\cla-\la I_n)\not=0$ and $\det(\cla^*-\la I_n)\not=0$, $w_\cla$ is invertible  and
\begin{align}& \label{D15}
w_\cla(x,\la)^{-1}= jw_\cla(x, \ov{\la})^*j.
\end{align}
Further we assume that
\begin{align}& \label{D16}
\det \cla\not=0,
\end{align}
and so $w_\cla(x,0)$ is well defined (in the points of invertibility of $\cls(x)$). 
We note that \eqref{D11} yields
\begin{align}& \label{D16+}
w_\cla^{\prime}(x,0)=-\wt q_0(x)w_\cla(x,0),
\end{align}
and we set
\begin{align}& \label{D17}
v(x,\la):=w_\cla(x,0)^{-1}w_\cla(x,\la).
\end{align}
Formulas \eqref{D11}, \eqref{D16+} and \eqref{D17}  imply that
\begin{align}& \label{D18}
v^{\prime}(x,\la)=\I \la j \wt H(x)v(x,\la)-\I \la v(x,\la)j  H(x),
\\ & \label{D19}
j\wt H(x)=w_\cla(x,0)^{-1}jH(x)w_\cla(x,0).
\end{align}
Thus, one can see that $j\wt H$ is linear similar to $jH$.
Moreover, in view of \eqref{D15} we can rewrite \eqref{D19} in the form
\begin{align}& \label{D20}
\wt H(x)=w_\cla(x,0)^{*}H(x)w_\cla(x,0).
\end{align}
According to \eqref{D20}, the equality $\wt H(x)=\wt H(x)^*$ is valid. Hence, $\wt H(x)$ is the {\it transformed} generalized Hamiltonian
of the {\it transformed} generalized canonical system
 \begin{align} &       \label{1.2'}
\wt w^{\prime}(x,\la)=\I \la j \wt  H(x)w(x,\la), \quad \wt H(x)=\wt H(x)^* \quad (x\geq 0).
 \end{align} 
Clearly, $\wt H \geq 0$ if $ H \geq  0$, and $\wt H >0$ if $ H >  0$.  Therefore, in the case of an initial canonical system,
the transformed system is also canonical.  By virtue of \eqref{D18}, a fundamental solution $\wt W$ of the transformed system
is given by the formula
\begin{align}& \label{D21}
\wt W(x,\la)=v(x,\la)W(x,\la),
\end{align}
where $W$ is a fundamental solution of the initial system.
\begin{Rk}\label{RkS} If $H(x)\geq 0$ and $\cls(0)>0$, the second equation in \eqref{D9} implies that
$\cls(x)>0$ for $x\geq 0$. In particular, $\cls(x)$ is invertible.
\end{Rk}
\begin{Rk}\label{RkD} In view of \eqref{D18} $($or \eqref{D21}$)$ the matrix function
$v(x,\la)$  is the so called Darboux matrix of the generalized canonical system.

According to \eqref{D14}, \eqref{D15} and \eqref{D17}, the representation of $v(x,\la)$
in terms of $\Lam(x)$ and $\cls(x)$ may be simplified. Namely, we have
\begin{align}\nn
v(x,\la)&=jw_\cla(x,0)^*jw_\cla(x,\la)
\\ & \label{D17+}
=I_m-\I\la j
 \Lam(x)^*(\cla^*)^{-1}\cls(x)^{-1} (\cla-\la I_n)^{-1}\Lam(x).
\end{align}
\end{Rk}
\section{Explicit solutions of the transformed \\  generalized canonical systems} \label{Tran}
\setcounter{equation}{0}
Consider the case, where the initial generalized Hamiltonian $H(x)$ has the form
\begin{align}& \label{E2}
H(x)=dj+\b(x)^*\b(x), \quad  \b(x):=\begin{bmatrix} \E^{\I c x}I_{m_1} & \E^{-\I c x}\a\end{bmatrix} .
\end{align}
Here, $\b(x)$ is an $m_1 \times m$ matrix function,  $\a$ is an $m_1 \times m_2$ matrix function and
\begin{align}& \label{E3}
c,\,d\in \BR; \quad \a\a^*=I_{m_1} \quad (m_2\geq m_1).
\end{align}
In view of \eqref{E2} and \eqref{E3}, we have
\begin{align}& \label{E3'}
\b(x)j\b(x)^*\equiv 0.
\end{align}
Recall that the matrix function $\Lam(x)$ is determined by $\Lam(0)$ and by the system
\begin{align}& \label{E1}
\Lam^{\prime}(x)=-\I \cla \Lam(x)j H(x).
\end{align} 
We construct {\it generalized eigenfunction} $\Lam(x)$ in the case \eqref{E2} explicitly.
\begin{Pn}\label{PnExpl} Let  \eqref{E2} and \eqref{E3} hold. Then, the matrix function $$\Lam(x)=\begin{bmatrix} \Phi_1(x) & \Phi_2(x) \end{bmatrix}$$
such that
\begin{align} \label{E4}
\Phi_1(x)=&\exp\{\I x(cI_n-d \cla)\}\big(\E^{\I xQ}f_1+\E^{-\I xQ}f_2\big),
\\  \nn
\Phi_2(x)=&\exp\{-\I x(cI_n+d \cla)\}\big(\E^{\I xQ}(\cla+cI_n+Q)\cla^{-1}f_1
\\ & \label{E5}
+\E^{-\I xQ}(\cla+cI_n-Q)\cla^{-1}f_2\big)\a,
\end{align}
where $f_k$ are $n\times m_1$ matrices, $Q$ is an $n\times n$ matrix and
\begin{align}& \label{E6}
\cla Q=Q\cla, \quad Q^2=c(2\cla+cI_n),
\end{align}
satisfies \eqref{E1}.
\end{Pn}
\begin{proof}. Using \eqref{E2}--\eqref{E5}
we derive
\begin{align}& \label{E7}
\Lam(x)j\b(x)^*=-\E^{-\I dx \cla}\big(\E^{\I xQ}(cI_n+Q)\cla^{-1}f_1
+\E^{-\I xQ}(cI_n-Q)\cla^{-1}f_2\big).
\end{align}
It follows from \eqref{E4} that
\begin{align}\nn
\frac{d}{d x}\Phi_1(x)=&-\I d \cla \Phi_1(x)+\I\exp\{\I x(cI_n-d \cla)\}
\\ & \label{E8}
\times \big(\E^{\I xQ}(cI_n+Q)f_1+\E^{-\I xQ}(cI_n-Q)f_2\big).
\end{align}
According to \eqref{E5}, we also have
\begin{align}\nn
\frac{d}{d x}\Phi_2(x)=&-\I d \cla \Phi_2(x)+\I\exp\{-\I x(cI_n+d \cla)\}
\\ & \nn \times
\big(\E^{\I xQ}(Q-cI_n)(\cla+cI_n+Q)\cla^{-1}f_1
\\ & \label{E9}
+\E^{-\I xQ}(Q+cI_n)(Q-\cla-cI_n)\cla^{-1}f_2\big).
\end{align}
Since $\Lam(x)=\begin{bmatrix} \Phi_1(x) & \Phi_2(x) \end{bmatrix}$ and $H(x)$ has the form \eqref{E2} (where \eqref{E3} holds)
relations \eqref{E6}--\eqref{E9} imply \eqref{E1}.
\end{proof}
It is easy to see that one can set $Q=0$ (in the Proposition \ref{PnExpl}) in the case $c=0$.
A more interesting case, where $c=0$ and \eqref{E6} holds, is generated by the matrices $\cla$ and $Q$ of the form
\begin{align}& \label{E9+}
\cla=\xi I_{2r}+\begin{bmatrix} 0 &\cla_{12} \\ 0 & 0\end{bmatrix} \quad (\xi\in \BC), \quad Q=\begin{bmatrix} 0 &Q_{12} \\ 0 & 0\end{bmatrix}
\end{align}
(where $\cla$ and $Q$ are $2r \times 2r$ matrices, $\cla_{12}$ and  $Q_{12}$ are some $r\times r$ matrices)
or by the block diagonal matrices with the blocks of the same form as the matrices on the right-hand sides of the equalities in \eqref{E9+}.

The next immediate corollary of \cite[Proposition B.1]{ALSgrav} (and its proof) deals with the case $c\not=0$.
\begin{Cy}\label{CyGrav} Let  $c\not=0$, let $\det(2\cla+cI_n)\not=0$, and let $\cle$ be the similarity transformation matrix and $\clj$
Jordan normal form in the representation
\begin{align}& \label{E10}
c(2\cla+cI_n)=\cle\clj \cle^{-1}.
\end{align}
Then, $Q$ satisfying \eqref{E6}
may be constructed explicitly and has the form
\begin{align}& \label{E11}
Q=\cle\cld \cle^{-1},
\end{align}
where $\cld$ is a block diagonal matrix with the blocks of the same orders as
the corresponding Jordan blocks of $\clj$. Moreover, the blocks of $\cld$ are upper
triangular Toeplitz matrices $($or scalars if the corresponding blocks of $\clj$ are
scalars$)$. If $z$ is the eigenvalue of some block of $\clj$, then the entries on the main diagonal
of the corresponding block of $\cld$ equal $\sqrt{z}$  $($and one can fix any of the
two possible values of $\sqrt{z}$ for this
main diagonal$)$.
\end{Cy}
Given generalized eigenfunction $\Lam(x)$, one can construct (explicitly)  the fundamental solution $\wt W(x,\la)$
of the transformed generalized canonical system
using relations \eqref{D13}, \eqref{D17}, \eqref{D21} and the second equality in \eqref{D9}. We note
that an explicit expression for $W(x,\la)$, which we need for this purpose, is constructed similar to the way
it is done in Proposition \ref{PnW}.
\section{The case of the spectral canonical systems} \label{Can}
\setcounter{equation}{0}
{\bf 1.} It follows from \eqref{D15}, \eqref{D20}, and \eqref{E3'} that the transformed generalized Hamiltonians constructed
in Section \ref{Tran} have the form
\begin{align}& \label{C1}
\wt H(x)= dw_\cla(x,0)^{*}jw_\cla(x,0)+\wt\b(x)^*\wt\b(x)=dj+\wt\b(x)^*\wt\b(x), \\
& \label{C2}
 \wt \b(x):=\b(x)w_\cla(x,0), \quad  \wt \b(x)j  \wt \b(x)^*=\b(x)j   \b(x)^*=0.
\end{align}
Here, $\wt \b$ is the corresponding transformation of $\b$. Setting 
\begin{align}& \label{C3}
d=0, \quad m_1=m_2=:p,
\end{align}
we obtain a class of canonical systems
\begin{align} &       \label{C4}
\wt  w^{\prime}(x,\la)=\I \la j \wt H(x)\wt w(x,\la), \quad \wt H(x)=\wt\b(x)^*\wt\b(x)\geq 0, \quad  \wt \b(x)j  \wt \b(x)^*=0.
 \end{align} 
 
Further in the text {\it we normalize the fundamental solutions $W$ and $\wt W$ of the  systems
\eqref{1.2} and  \eqref{1.2'}, respectively, setting}
\begin{align}& \label{C6}
W(0,\la)=\wt W(0,\la)=I_{2p}.
\end{align}
We write down the Hamiltonian $H(x)$ given by \eqref{E2}, \eqref{E3}, and \eqref{C3} in the form
\begin{align}& \label{C7}
H(x)=\E^{-\I cx j}\clk\E^{\I cx j}, \quad \clk:=\begin{bmatrix} I_p & \a \\ \a^* & I_p \end{bmatrix} \quad ( \a\a^*=I_p).
\end{align}
\begin{Pn}\label{PnW}
The fundamental solution of the canonical system  \eqref{1.3}, where $\Im(\la)\not=0$, the Hamiltonian $H$ is given by \eqref{C7} and $c\not=0$, has the form
\begin{align}& \label{C8}
W(x,\la)=\E^{-\I cxj}E(\la)\begin{bmatrix}\E^{\I z_1(\la) x} I_p & 0\\ 0 & \E^{\I z_2(\la) x} I_p \end{bmatrix}E(\la)^{-1}, \quad E=\begin{bmatrix}E_1 & E_2 \end{bmatrix},
\\
& \label{C9}
 E_i:=\begin{bmatrix}-\a \\ \frac{1}{\la}(\la+c-z_i)I_p \end{bmatrix}, \quad z_{i}^2=c(2\la+c) \quad (i=1,2), \quad \Im(z_1)>0.
\end{align}
\end{Pn}
\begin{proof}. It is easy to see that $E$ is invertible (one may consider, for instance, the linear span of the rows of $E$, which coincides
with $\BC^{2p}$). Moreover, using the equality 
$$\la -\frac{1}{\la}(\la+c)(\la+c-z_i)=\frac{z_i}{\la}(\la+c-z_i),$$ 
we have 
\begin{align}& \label{C10}
EZE^{-1}E=EZ=(\la j\clk+cj)E \quad {\mathrm{for}} \quad Z=\diag\{z_1I_p, \, z_2 I_p\}. 
\end{align}
It follows that
\begin{align}& \label{C11}
EZE^{-1}=\la j\clk+cj.
\end{align}
Relations \eqref{C8} and \eqref{C11} yield
\begin{align}& \label{C12}
W(x,\la)=\E^{-\I cxj}\E^{\I  x (\la j\clk+cj)},
\end{align}
and for $W(x,\la)$ of the form \eqref{C12} we immediately obtain
\begin{align}& \label{C13}
W^{\prime}(x,\la)=\I \la j\E^{-\I cxj}\clk\E^{\I cxj}W(x,\la), \quad W(0,\la)=I_{2p}.
\end{align}
Taking into account \eqref{C7} and \eqref{C13}, we see that $W$ given by \eqref{C8}, \eqref{C9} is, indeed, the normalized
fundamental solution of the canonical system described in the proposition.
\end{proof}
{\bf 2.} Further in this section, we assume that
\begin{align}& \label{C13+}
\cls(0)>0, \quad c\not= 0,
\end{align}
so that the statements of Remark \ref{RkS} and Proposition \ref{PnW} may be used.
After normalization \eqref{C6} formula \eqref{D21} takes the form
\begin{align}& \label{C14}
\wt W(x,\la)=v(x,\la)W(x,\la)v(0,\la)^{-1}.
\end{align}
For system \eqref{C1}--\eqref{C3}, in view of \eqref{D17}, \eqref{C8} and \eqref{C14} we obtain
 \begin{align}& \label{C15}
\wt \b(x)\wt W(x,\la)v(0,\la)E_1(\la)=\E^{\I z_1(\la)x}\b(x)w_\cla(x,\la)\E^{-\I c x j}E_1(\la).
\end{align}
Taking into account \eqref{C15} (and some definitions and considerations on Weyl--Titchmarsh theory in Appendix \ref{Weyl}), we derive the following theorem.
\begin{Tm} \label{TmPhi} Canonical system with Hamiltonian of the form \eqref{C1}--\eqref{C3} on $[0, \, \infty)$, where \eqref{C13+} holds, has
a unique Weyl function $($Weyl's limit point case$)$. This Weyl function is given explicitly by the
formula
 \begin{align}& \label{C16}
\vp(\la)=\begin{bmatrix}0 & I_p\end{bmatrix}v(0,\la)E_1(\la)\big(\begin{bmatrix}I_p & 0\end{bmatrix}v(0,\la)E_1(\la)\big)^{-1},
\end{align}
where $E_1$ has the form \eqref{C9}.
\end{Tm}
\begin{proof}. First, we note that formulas \eqref{D14}--\eqref{D17} (and \cite[Corollary E.3]{SaSaR}) yield
\begin{align}& \label{C17}
v(0,\la)^*jv(0,\la)\geq j, \quad  v(0,\la)jv(0,\la)^*\geq j \quad (\la \in \BC_+).
\end{align}
It easily follows from \eqref{C9} and \eqref{C17} that
\begin{align}& \label{C18}
E_1\Big(-\frac{c}{2}+\ve \I\Big)^*j E_1\Big(-\frac{c}{2}+\ve \I\Big)>0,
\quad \begin{bmatrix}I_p & 0\end{bmatrix}v(0,\la)jv(0,\la)^*\begin{bmatrix}I_p \\ 0\end{bmatrix}>0, 
\end{align}
where the first inequality holds (at least) for small $\ve>0$ and the second inequality holds for
all $\la \in \BC_+$ (excluding the part of spectrum of $\cla$ situated in $\BC_+$). Hence (see, e.g., \cite[Proposition 1.43]{SaSaR}), $\det \big(\begin{bmatrix}I_p & 0\end{bmatrix}v(0,\la)E_1(\la)\big)\not\equiv 0$,
and so 
\begin{align}& \label{C19}
\det \big(\begin{bmatrix}I_p & 0\end{bmatrix}v(0,\la)E_1(\la)\big)\not= 0 \quad {\mathrm{for}} \quad \la \in \BC_+,
\end{align}
excluding, possibly, some isolated points. In other words, $\vp(\la)$ in \eqref{C16} is well defined.

Now, we will show that for such $\la$ that $\Im\big(z_1(\la)\big)$ is sufficiently large (excluding, possibly,
isolated points) the relation
\begin{align}& \label{C20}
\wt \b(x)\wt W(x,\la)\begin{bmatrix}I_p \\ \vp(\la)\end{bmatrix}\in L_2^{p\times p}(0,\infty)
\end{align}
is valid. Indeed, the matrix functions $\b(x)$ and $\E^{-\I c x j}$ on the right-hand side of  \eqref{C15}
are bounded. In view of  \eqref{D9}, we have
$$\int_0^r \cls(x)^{-1}\Lam(x)j\b(x)^*\b(x)j\Lam(x)^*\cls(x)^{-1}=\cls(0)^{-1}-\cls(r)^{-1}\leq \cls(0)^{-1}.$$
Therefore, we obtain
\begin{align}& \label{C21}
 \b(x)j\Lam(x)^*\cls(x)^{-1} \in L_2^{p\times n}(0,\infty).
\end{align}
Finally, Proposition \ref{PnExpl} and Corollary \ref{CyGrav} show that the matrix function
$\E^{\I z_1(\la)x}\Lam(x)$ is bounded for sufficiently large values of $\Im\big(z_1(\la)\big)$.
Taking into account the definition \eqref{D13} of $w_\cla$ and considerations above,
we see that the right-hand side of \eqref{C15} belongs $ L_2^{p\times p}(0,\infty)$
(for sufficiently large values of $\Im\big(z_1(\la)\big)$).  Thus, the  left-hand side of \eqref{C15} belongs $ L_2^{p\times p}(0,\infty)$
as well, and so
\eqref{C20} holds
for $\vp(\la)$ given by \eqref{C16}.

Assume that for some $\la=\la_0 \in \BC_+$ we have \eqref{C20} and also have
\begin{align}& \label{C22}
\wt \b(x)\wt W(x,\la_0)\begin{bmatrix}I_p \\ \wh \vp(\la_0)\end{bmatrix}\in L_2^{p\times p}(0,\infty), \quad  \vp(\la_0)\not=\wh \vp(\la_0).
\end{align}
We will show (by contradiction) that this is impossible for sufficiently large values of $\Im\big(z_1(\la_0)\big)$.
Indeed, since \eqref{C20} implies 
$$\wt \b(x)\wt W(x,\la_0)v(0,\la_0)E_1(\la_0) \in L_2^{p\times p}(0,\infty),$$
additional relations \eqref{C22}  yield the existence of $f\in \BC^p$ such that
\begin{align}& \label{C23}
\wt \b(x)\wt W(x,\la_0)v(0,\la_0)E_2(\la_0)f \in L_2^{p\times 1}(0,\infty) \quad (f\not=0).
\end{align}
On the other hand, taking into account that $z_2(\la)=-z_1(\la)$ we similar to \eqref{C15} derive
 \begin{align}& \label{C24}
\wt \b(x)\wt W(x,\la)v(0,\la)E_2(\la)=\E^{-\I z_1(\la)x}\b(x)w_\cla(x,\la)\E^{-\I c x j}E_2(\la).
\end{align}
Next, we should consider $g(x,\la)=\b(x)w_\cla(x,\la)\E^{-\I c x j}$ in a more detailed way, and we note that
according to \eqref{D9}, \eqref{D13}, \eqref{E2}, \eqref{E3}, \eqref{E4}, and \eqref{E5} the entries $g_{ik}$ of $g$
 admit representation
 \begin{align}& \label{C25}
g_{ik}(x,\la)=\sum_{s=1}^{N_1} P_s(\la)x^{\ell_s}\E^{h_s x}\big/\left(P(\la)\sum_{s=1}^{N_2}x^{n_s}\E^{\zeta_s x}\right).
\end{align}
where $P$ and $P_s$ are polynomials, and $N_1$, $P_s$, $\ell_s$ and $h_s$ depend on $i,\, k$.
Moreover,  similar to \eqref{C19} one can show that (excluding isolated points $\la$) we have
 \begin{align}& \nn
\det\big(g(x,\la)\wh E(\la)\big)\not=0, \quad \wh E(\la):=\begin{bmatrix}-\a \\ \big((\la+c)/{\la} \big)I_p \end{bmatrix},
\end{align}
where $\wh E$ is the ``rational part" of $E_2$. It follows that
\begin{align}& \label{C26}
g(x,\la)E_2(\la)f\not=0.
\end{align}
Taking into account \eqref{C24}--\eqref{C26}, we see that \eqref{C23} (and so \eqref{C22}) does not hold for sufficiently
large values of $\Im\big(z_1(\la_0)\big)$.

Since \eqref{C22} does not hold for sufficiently
large values of $\Im\big(z_1(\la_0)\big)$ (excluding, may be, isolated points), there is an open domain in $\BC_+$,
where $\vp(\la)$ (given by \eqref{C16}) is uniquely defined via \eqref{C20}. Thus, 
each Weyl function of our system coincides with $\vp(\la)$ in this domain (see the Definition \ref{DnW} of the Weyl  functions). 
Recall that Weyl functions are holomorphic
in $\BC_+$. Hence, the Weyl function of our system is unique (and its existence follows from Proposition \ref{PnH}).
We see that the Weyl function exists, is unique and coincides with $\vp(\la)$ in some domain.
Therefore, $\vp(\la)$ given by  \eqref{C16} is the Weyl function and admits holomorphic continuation in all
$\BC_+$.
\end{proof}

\begin{Rk}\label{RkC0} In Proposition \ref{PnW} and Theorem \ref{TmPhi}, we assume that $c\not=0$.
The constructions are much simpler when $c=0$. In particular, we recall that $\b j\b^*\equiv 0$ $($see 
\eqref{E3'}$)$. Moreover, $\clk$ given  in \eqref{C7} equals $\b^*\b$ as $c=0$.
Hence, $\clk j\clk=0$. Therefore, using \eqref{C7} we have:
$$H(x)\equiv \clk, \quad W(x,\la)=\E^{\I \la x j \clk}= I_{2p}+\I \la x j \clk$$
for the case $c=0$.
\end{Rk}
\section{Matrix string equation} \label{MSchr}
\setcounter{equation}{0}
Consider again the case of the initial canonical systems 
$$ w^{\prime}(x,\la)=\I \la j H(x)w(x,\la),$$ 
where $ H(x)= \b(x)^*\b(x)$ and $\b(x)$
are $p \times 2p$ matrix functions.  According to \eqref{D20},  the transformed Hamiltonians
(of the GBDT-transformed canonical systems  \eqref{1.2'}) have the form
\begin{align}& \label{S1}
\wt H(x)=\wt \b(x)^*\wt \b(x), \quad \wt \b(x)=\b(x)w_\cla(x,0).
\end{align}
When the matrix functions $\b(x)$ have the form presented in \eqref{E2} (and \eqref{E3}, \eqref{C3} hold),
our assertions below show (in view of  Appendix \ref{String}) that the considered  transformed
canonical systems correspond to a special subclass of string equations. Thus, our explicit
formulas may be transferred for the case of string equations as explained in Remark \ref{RkFS}.
\begin{Pn}\label{PnS}
Let $\b(x)$ satisfy the equality
\begin{align}& \label{S2}
 \b(x)j\b(x)^*=0.
\end{align}
Then, $\wt \b(x)$ satisfies the relations
\begin{align}& \label{S3}
\wt  \b(x)j \wt \b(x)^*=0, \quad \wt \b(x)^{\prime}j \wt\b(x)^*=\b(x)^{\prime}j\b(x)^*.
\end{align}
\end{Pn}
\begin{proof}.  Recall (see, e.g., \eqref{D15}) that
\begin{align}& \label{S3'}
w_\cla(x,0)jw_\cla(x,0)^*=j.
\end{align}
The first equality in \eqref{S3} easily follows  from \eqref{S3'} (and was already stated in \eqref{C2}).
Formulas  \eqref{D16+}, \eqref{S3'} and the second equality
in \eqref{S1} imply that
\begin{align}\nn
 \wt \b(x)^{\prime}j \wt\b(x)^*&=\b(x)^{\prime}w_\cla(x,0)jw_\cla(x,0)^*\b(x)^*+\b(x)w_\cla^{\prime}(x,0)jw_\cla(x,0)^*\b(x)^*
 \\  \label{S4} &
=\b(x)^{\prime}j\b(x)^*-\b(x)\wt q_0(x)j\b(x)^*.
\end{align}
The definition \eqref{D12} of $\wt q_0$ and the  equality  \eqref{S2} yield
\begin{align}& \label{S5}
\b(x)\wt q_0(x)j\b(x)^*=0.
\end{align}
The second equality in  \eqref{S3}  is immediate from \eqref{S4} and \eqref{S5}.
\end{proof}
\begin{Cy}\label{CyS} Let $\b(x)$ be given by \eqref{E2}, where $c=\frac{1}{2}$ and $\a\a^*=I_p$. Then, \eqref{S2} holds and
\begin{align}& \label{S6}
 \wt \b(x)^{\prime}j \wt\b(x)^*=\I I_p.
\end{align}
\end{Cy}

\section{Examples and applications} \label{ExAp}
\setcounter{equation}{0}
{\bf 1.} Let us consider explicit examples of the Hamiltonians $\wt H(x)=\wt \b(x)^*\wt\b(x)$, corresponding
Darboux matrices $v(x,\la)$, fundamental solutions $\wt W(x,\la)$, and Weyl functions $\vp(\la)$.
\begin{Ee} \label{Ee1} In our first example, we assume that
\begin{align}& \label{P1}
p=n=1, \quad \cla=a\not=\ov{a} \quad (a \in \BC), \quad c\not=0, \quad d=0,
\end{align}
where the condition $a\not=\ov{a}$ provides an easy recovery  of $\cls(x)$ from \eqref{D10}.
\end{Ee}
Recall that according to the second equalities in \eqref{E2}, \eqref{E3}, and \eqref{E6},  we have
\begin{align}& \label{P2}
\b(x):=\begin{bmatrix} \E^{\I c x} & \E^{-\I c x}\a\end{bmatrix}, \quad |\a|=1, \quad Q=\sqrt{2ac+c^2}.
\end{align}
In order to define the sign of the square root above, we assume that $\Im(Q)>0$. 
By virtue of \eqref{E4}, \eqref{E5}) and \eqref{P1}, we obtain
\begin{align}\nn
a\Lam(x)=&\begin{bmatrix}a\big({f_1}\E^{\I x {Q}}+{f_2}\E^{-\I x {Q}}\big) &
\a \big({(a+c+Q)f_1}\E^{\I x {Q}}+{(a+c-Q)f_2}\E^{-\I x {Q}}\big)
\end{bmatrix}
\\ & \label{P4+}\times 
\E^{\I c x j},
\end{align}
where $f_1$ and $f_2$ are scalars and $a\Lam$ is written down more conveniently than $\Lam$.
It follows from \eqref{D10} and \eqref{P4+} that
\begin{align}\nn
\cls(x)=&\frac{\I}{a-\ov{a}}\Big(\big|f_1\E^{\I x Q}+f_2\E^{-\I x Q}\big|^2
\\ & \label{P3}
-\frac{1}{|a|^2}\big|(a+c+Q)f_1\E^{\I x Q}+(a+c-Q)f_2\E^{-\I x Q}\big|^2
\Big),
\end{align}
and the requirement $\cls(0)>0$ takes the form
\begin{align}& \label{P4}
\I(\ov{a}-a)\big(|a(f_1+f_2)|^2-\big|(a+c+Q)f_1+(a+c-Q)f_2\big|^2\big)>0.
\end{align}
Relations \eqref{D13},  \eqref{C2} and \eqref{P1}--\eqref{P4+} yield
\begin{align}\nn
\wt \b(x)=&\b(x)-\frac{\I}{a|a|^2\cls(x)}\Big(\ov{a}\big(\ov{f_1}\E^{-\I x \ov{Q}}+\ov{f_2}\E^{\I x \ov{Q}}\big)
\\  \label{P5} & 
- \ov{\a}\big(\ov{(a+c+Q)f_1}\E^{-\I x \ov{Q}}+\ov{(a+c-Q)f_2}\E^{\I x \ov{Q}}\big)\Big)\big(a\Lam(x)\big).
\end{align}
According to \eqref{D17+} and \eqref{P1}, the corresponding Darboux matrix is given by the formula
\begin{align}&\label{P6}
v(x,\la)= I_2-\frac{\I \la}{\ov{a}|a|^2(a-\la)\cls(x)}j\big(a\Lam(x)\big)^*\big(a\Lam(x)\big).
\end{align}
Formulas \eqref{C8}, \eqref{C14} and \eqref{P4+}, \eqref{P6} give explicitly
fundamental solutions of the canonical systems with $\wt \b$ of the form \eqref{P5}.
In view of \eqref{C16} and \eqref{P6}, the Weyl functions $\vp$ of such canonical
systems on $[0,\,\infty)$ have the form:
\begin{align} \label{P7}
\vp(\la)=&{\psi_1(\la)}\big/{\psi_2(\la)}, \\
\nn
 \psi_1(\la)=&\ov{a}|a|^2\cls(0)(a-\la)\big(\la+c-z_1(\la)\big)
\\ & \label{P8} 
 +\I\ov{\a}\big((\ov{a}+c+\ov{Q})\ov{f_1}+(\ov{a}+c-\ov{Q})\ov{f_2}\big)\la
h(\la),
\\  \label{P9} 
\psi_2(\la)=&\a \ov{a}|a|^2\cls(0)(\la-a)\la-\I \ov{a}(\ov{f_1}+\ov{f_2})\la h(\la),
\end{align}
where $z_1(\la)=\sqrt{c(2\la+c)} \,\, (\Im(z_1)>0)$,
\begin{align}
 \nn
h(\la):=a\Lam(0)E_1(\la)=& \a\big({(a+c+Q)f_1}+{(a+c-Q)f_2}\big)\big(\la+c-z_1(\la)\big)
\\ &  \label{P10} 
-\a a(f_1+f_2)\la .
\end{align}
\begin{Ee} \label{Ee2} Now, assume that 
\begin{align}& \label{P11}
p=1, \quad n=2, \quad c=0, \quad d=0, \quad \cla=\begin{bmatrix} \xi & a\\ 0 & \xi \end{bmatrix} \quad (\xi\in \BR, \,\,\xi\not= 0), 
\\
& \label{P11+} 
Q=\begin{bmatrix} 0 & q \\ 0 & 0\end{bmatrix}, \quad f_1=\begin{bmatrix} f \\  0\end{bmatrix} \quad f_2=  \begin{bmatrix} 0 \\  g\end{bmatrix};
\quad
q,f,g\in \BC,
 \quad f\not=0,
\quad g\not=0.
\end{align}
\end{Ee}
In this case, we have
\begin{align}& \label{P12}
\b=\begin{bmatrix} 1 &  \a\end{bmatrix}, \quad \E^{\pm \I x Q}=I_2\pm \I x Q, \\
& \label{P13}
 (\cla-\la I_2)^{-1}=(\xi-\la)^{-1}I_2-(\xi-\la)^{-2}\begin{bmatrix} 0 & a\\ 0 & 0\end{bmatrix} .
\end{align}
Hence, formulas \eqref{E4}, \eqref{E5} and simple calculations yield
\begin{align}& \label{P14}
\Lam(x)=\begin{bmatrix} f-\I qgx & -qg\big(\I x +\xi^{-1}\big) \\ g & g\end{bmatrix}
\end{align}
In view of \eqref{P14}, the required matrix identity \eqref{D8}  may be written in the form
\begin{align}& \label{P15}
\begin{bmatrix} a\cls_{21}(0)-\ov{a}\cls_{12}(0) & a\cls_{22}(0)\\ -\ov{a}\cls_{22}(0) & 0\end{bmatrix}=\I
\begin{bmatrix} |f|^2-|qg\xi^{-1}|^2 & \ov{g}\big(f +gq\xi^{-1}\big) \\ g \big(\ov{f} +\ov{gq}\xi^{-1}\big)& 0\end{bmatrix},
\end{align}
where $\cls_{ik}$ are the entries of $\cls$.  Hence, we cannot choose an arbitrary entry $a$ in $\cla$ but demand $f +gq\xi^{-1}\not=0$
and choose $a$ and $\cls_{22}(0)$ satisfying the following conditions (which is always possible):
\begin{align}& \label{P16}
a\cls_{22}(0)=\ov{g}\big(f +gq\xi^{-1}\big),\quad a\not=0, \quad \cls_{22}(0)>0.
\end{align}
Next, we choose $\cls_{12}(0)$ (and so $\cls_{21}(0)=\ov{\cls_{12}(0)}$) such that \eqref{P15} holds, and we choose such $\cls_{11}(0)>0$
that $\cls(0)>0$.  

Since $\xi \in \BR$, we cannot use \eqref{D10} in order to recover $\cls(x)$ from $\Lam(x)$
and construct $\cls(x)$ in a different way.
It follows from \eqref{E2}, \eqref{P11} and \eqref{P14} that
\begin{align}& \label{P17}
\Lam(x)j\b^*=\begin{bmatrix} C_1x+C_2 \\ C_3\end{bmatrix}, \quad C_1=\I(\ov{\a}-1)qg, \quad C_2=f+\ov{\a}qg\xi^{-1}, \\
& \label{P18}
 C_3=g(1-\ov{\a}).
\end{align}
Therefore, the second equality in \eqref{D9} yields
\begin{align}\label{P19}
\cls(x)&=\cls(0)+\int_0^x \Lam(t)j\b^*\big(\Lam(t)j\b^*\big)^*dt
\\ & \nn
=\cls(0)+\begin{bmatrix} \frac{1}{3}|C_1|^2x^3+\Re\big(C_1\ov{C_2}\big)x^2+|C_2|^2x
& \frac{1}{2}C_1\ov{C_3}x^2+C_2\ov{C_3}x
\\ \frac{1}{2}\ov{C_1}{C_3}x^2+\ov{C_2}{C_3}x & |C_3|^2x
\end{bmatrix}.
\end{align}
Using \eqref{P17},  we rewrite the equality  \eqref{C2} for $\wt \b$ (transformed $\b$) in the form
\begin{align}& \label{P20}
\wt \b(x)=\begin{bmatrix} 1 & \a\end{bmatrix}-\I \begin{bmatrix}\, \ov{C_1}x+\ov{C_2} &\quad  \ov{C_3}\,\end{bmatrix}\cls(x)^{-1}\cla^{-1}\Lam(x),
\end{align}
where $\cls(x)$, $\cla$ and $\Lam(x)$ are given in  \eqref{P19}, \eqref{P11} and \eqref{P14}, respectively.
Finally, the Darboux matrix $v(x,\la)$ is expressed via $\Lam(x)$ and $\cls(x)$ in \eqref{D17+},
and the expression for the corresponding fundamental solution $\wt W$ follows from
\eqref{C14} and Remark \ref{RkC0}.

{\bf 2.} Relations \eqref{D9} and \eqref{D10} imply an important equality (see \cite[(2.13)]{ALS17}):
\begin{align}& \label{P21}
\big(\Lam^*\cls^{-1}\big)^{\prime}=\I Hj\Lam^*\cls^{-1}\cla+\wt q_0^{\, *}\Lam^*\cls^{-1}.
\end{align}
We assume that $\cls(0)>0$ and $H(x)\geq 0$, that is, $\cls(x)>0$  for $x\geq 0$,
and so $\cls(x)^{-1}$ is well defined (see Remark \ref{RkS}).
In view of \eqref{D15}, \eqref{D16+} and \eqref{P21}, for $\wt H$ of the form \eqref{D20} and $Y$ given by
\begin{align}& \label{P22}
Y(x,t)=jw_\cla(x,0)^*\Lam(x)^*\cls(x)^{-1}\E^{\I t \cla},
\end{align}
we have
\begin{align}& \label{P23}
\wt H(x) \frac{\p}{\p t}Y(x,t)=j\frac{\p}{\p x}Y(x,t) \quad (m_1=m_2=p), \quad  x\geq 0.
\end{align}
In other words, the $2p \times n$ matrix function $Y$ (or, equivalently, the columns of $Y$) satisfies the dynamical
canonical system \eqref{P23}.

Taking into account \eqref{D10} and \eqref{D13}, we rewrite $w_\cla(x,0)^*\Lam(x)^*\cls(x)^{-1}$ 
in a simpler form (in terms of $\Lam(x)$ and $\cls(x)$):
\begin{align}& \nn
w_\cla(x,0)^*\Lam(x)^*\cls(x)^{-1}=\Lam(x)^*\big(\cla^{*}\big)^{-1}\cls(x)^{-1}\cla^{-1}.
\end{align}
Hence, 
\begin{align}& \label{P24}
Y(x,t)=j\Lam(x)^*\big(\cla^{*}\big)^{-1}\cls(x)^{-1}\E^{\I t \cla}\cla^{-1}.
\end{align}
\begin{Pn}\label{DCAnS} Let the initial Hamiltonian $H(x)\geq 0$ be given, and let the relations \eqref{D8},
$\cls(0)>0,$ and $\det  \cla \not=0$ hold. Then, $Y$ of the form \eqref{P24} satisfies 
dynamical canonical system \eqref{P23}, where the the transformed Hamiltonian $\wt H$
is given by \eqref{D20}.
\end{Pn}

In this way, explicit expressions for $\Lam$ and $\cls$ in Examples \ref{Ee1} and \ref{Ee2}
give us explicit expessions for $Y(x,t)$. Moreover, it is immediate from \eqref{P11}
that $\E^{\I t \cla}$ in \eqref{P24} takes under assumptions of  Example \ref{Ee2}  a simple form
\begin{align}& \label{P25}
\E^{\I t \cla}=\E^{\I t \xi}\left(I_2+\I t a\begin{bmatrix} 0 &  1\\ 0 & 0\end{bmatrix}\right).
\end{align}

\appendix
\section{Canonical systems: \\ Weyl--Titchmarsh theory} \label{Weyl}
\setcounter{equation}{0}
Consider generalized canonical system \eqref{1.2}. It is immediate  that
the fundamental solution $W$
of \eqref{1.2}  satisfies the equality
\begin{align}& \label{W0}
\frac{d}{dx}\big(W(x,\ov{\mu})^*jW(x,\la)\big)=\I({\la - \mu})W(x,\la)^*H(x)W(x,\la).
\end{align}
In view of  \eqref{W0} (for the case $\mu=\ov{\la}$) and of the normalization $W(0,\la)=I_m$ always assumed in this appendix, we have
\begin{align}& \label{W1}
\int_0^rW(x,\la)^*H(x)W(x,\la)dx=\frac{\I}{\la - \ov{\la}}\big(j-W(r,\la)^*jW(r,\la)\big),
\end{align}
for $\la\not\in \BR$ and $r\geq 0$. 
Moreover,  \eqref{W0} for the case $\mu={\la}$ implies that
\begin{align}& \label{W2}
W(r,\ov{\la})^*jW(r,\la)\equiv  j\equiv W(r,\la)jW(r,\ov{\la})^*.
\end{align}
Further in the appendix, we will deal with the general-type (i.e., not necessarily related to explicit solutions) canonical  system
 \eqref{1.3} on $[0,\, \infty)$.
Since $H\geq 0$, formula \eqref{W1}  yields
\begin{align}& \label{W3}
W(r_2,{\la})^*jW(r_2,\la)\leq W(r_1,{\la})^*jW(r_1,\la)\leq j,
\\ & \label{W4}
j\leq W(r_1,\ov{\la})^*jW(r_1,\ov{\la})\leq W(r_2,\ov{\la})^*jW(r_2,\ov{\la}) \quad (r_1 \leq r_2, \quad \la \in \BC_+). 
\end{align}
Next, introduce the families $\cln(r)$ of linear-fractional (M\"obius) transformations
\begin{align} \nn
\phi(r,\la)=&\big(\clw_{21}(r,\la)\clp_1(\la)+\clw_{22}(r,\la)\clp_2(\la)\big)
\\ & \label{W5}
\times\big(\clw_{11}(r,\la)\clp_1(\la)+\clw_{12}(r,\la)\clp_2(\la)\big)^{-1},
\end{align}
where $\clw_{ik}$ and $\clp_k$ are $p\times p$ matrix functions,
\begin{align} & \label{W6}
\clw(r,\la)=\{\clw_{ik}(r,\la)\}_{i,k=1}^2:=jW(r,\ov{\la})^*j,
\end{align}
and $\clp_1(\la)$, $\clp_2(\la)$  are pairs of meromorphic in $\BC_+$ matrix functions
(so called nonsingular pairs with property-$j$) such that
\begin{align} & \label{W7}
\clp_1(\la)^*\clp_1(\la)+\clp_1(\la)^*\clp_1(\la)>0, \quad  \begin{bmatrix}\clp_1(\la)^* & \clp_2(\la)^*\end{bmatrix}j \begin{bmatrix}\clp_1(\la) \\ \clp_2(\la)\end{bmatrix}\geq 0,
\end{align}
where the first inequality holds in one point (at least) of $\BC_+$, and the second inequality holds in all the points of analyticity
of $\clp_k$ $(k=1,2)$. It follows from \eqref{W7} by contradiction that
\begin{align} & \label{W8}
\det\big(\clw_{11}(r,\la)\clp_1(\la)+\clw_{12}(r,\la)\clp_2(\la)\big)\not\equiv 0.
\end{align}
Indeed, formulas \eqref{W2}, \eqref{W3} and \eqref{W6} imply that $\clw(r,{\la})^*j \clw(r,\la)\geq j$, which yields
\begin{align}& \label{Z1}
\begin{bmatrix}\clp_1(\la)^* & \clp_2(\la)^*\end{bmatrix}\clw(r,{\la})^*j \clw(r,\la)\begin{bmatrix}\clp_1(\la) \\ \clp_2(\la)\end{bmatrix}\geq 0
\end{align}
for such $\la$ in $\BC_+$  that \eqref{W7} holds.  On the other hand, if we have 
$$\det\big(\clw_{11}(r,\la)\clp_1(\la)+\clw_{12}(r,\la)\clp_2(\la)\big)= 0,$$
 then (for some $g\not=0$)
\begin{align} & \label{Z2}
\big(\clw_{11}(r,\la)\clp_1(\la)+\clw_{12}(r,\la)\clp_2(\la)\big)g= 0 \quad (g\in \BC^p),
\end{align}
and so we obtain
\begin{align} & \label{Z3}
g^*\begin{bmatrix}\clp_1(\la)^* & \clp_2(\la)^*\end{bmatrix}\clw(r,{\la})^*j \clw(r,\la)\begin{bmatrix}\clp_1(\la) \\ \clp_2(\la)\end{bmatrix}g<0.
\end{align}
Clearly, \eqref{Z3} contradicts \eqref{Z1}.

Let us rewrite \eqref{W5} in the form 
$$\begin{bmatrix}I_p \\ \phi(\la)\end{bmatrix}=jW(r,\ov{\la})^*j\begin{bmatrix}\clp_1(\la) \\ \clp_2(\la)\end{bmatrix}
\big(\clw_{11}(r,\la)\clp_1(\la)+\clw_{12}(r,\la)\clp_2(\la)\big)^{-1}.$$
Now, setting 
\begin{align} & \label{W9}
\mathfrak{A}(r,\la):=W(r,\la)^*jW(r,\la),
\end{align}
and using \eqref{W2}, we see that formulas \eqref{W5}--\eqref{W7} (i.e., the relation $\phi\in \cln(r)$) yield
\begin{align} & \label{W10}
\begin{bmatrix}I_p & \phi(\la)^*\end{bmatrix}\mathfrak{A}(r,\la)\begin{bmatrix}I_p \\ \phi(\la)\end{bmatrix}\geq 0.
\end{align}
Moreover, according to \eqref{W3}, \eqref{W9} and \eqref{W10}, $\phi(\la)$ is holomorphic and contractive in $\BC_+$. On the other hand, if $\phi$ satisfies \eqref{W10}, we set
\begin{align} & \label{W11}
\begin{bmatrix}\clp_1(\la) \\ \clp_2(\la)\end{bmatrix}=W(r,\la)\begin{bmatrix}I_p \\ \phi(\la)\end{bmatrix},
\end{align}
and see that the relations \eqref{W5}--\eqref{W7} are valid. Thus, 
\begin{align} & \label{W12}
\phi(\la)\in \cln(r)
\end{align}
is equivalent to \eqref{W10}.  Therefore, according to \eqref{W3}, $\cln(r_2)$ is embedded in $\cln(r_1)$:
\begin{align} & \label{W13}
\cln(r_2) \subseteq \cln(r_1) \quad (r_1<r_2).
\end{align}
By virtue of Montel's theorem, there is a sequence $\{\phi_k(\la)\}$  such that
\begin{align} & \label{W14-}
\phi_k\in \cln(r_k), \quad r_k\to \infty \quad ({\mathrm{for}} \,\, k\to \infty) ,
\end{align}
and $\phi_k(\la)$ tend uniformly (on any compact in $\BC_+$) to some matrix function $\vp(\la)$.
Thus, $\vp(\la)$ is holomorphic and satisfies \eqref{W10} for any $r>0$. In other words,
\begin{align} & \label{W14}
\vp(\la)\in \bigcap_{r>0}\cln(r).
\end{align}

Let us write down $\cln(r)$ in the Weyl matrix disk form. Taking into account \eqref{W3}
and \eqref{W9}, we obtain
\begin{align} & \label{Z4}
-\mf{A}_{22}(r_2,\la)\geq -\mf{A}_{22}(r_1,\la)\geq I_p \quad (r_2>r_1); 
\\ & \label{Z5-}
\mf{A}(r_2,\la)^{-1}\geq \mf{A}(r_1,\la)^{-1}\geq j,
\\ & \label{Z5}
\big(\mf{A}(r,\la)^{-1}\big)_{11}=\big(\mf{A}_{11}(r,\la)-\mf{A}_{12}(r,\la)\mf{A}_{22}(r,\la)^{-1}\mf{A}_{21}(r,\la)\big)^{-1}\geq I_p,
\end{align}
where $\mf{A}_{ik}(r,\la)$ and $\big(\mf{A}(r,\la)^{-1}\big)_{ik}$, respectively, are $p\times p$ blocks of $\mf{A}(r,\la)$ and $\mf{A}(r,\la)^{-1}$.
The invertibility of  $\mf{A}_{11}-\mf{A}_{12}\mf{A}_{22}^{-1}\mf{A}_{21}$ in \eqref{Z5} follows from the invertibility of
$\mf{A}$ and $\mf{A}_{22}$ (for this and for the equality in \eqref{Z5} see \cite[p. 21]{SaL2}).
In particular, we derive from \eqref{Z4} and \eqref{Z5} that the following positive definite square roots are uniquely
defined:
\begin{align} & \label{Z6}
\rho_L(r,\la)=\big(-\mf{A}_{22}(r_2,\la)^{-1}\big)^{1/2},
\\ & \label{Z7}
 \rho_R(r,\la)=\big(\mf{A}_{11}(r,\la)-\mf{A}_{12}(r,\la)\mf{A}_{22}(r,\la)^{-1}\mf{A}_{21}(r,\la)\big)^{1/2}.
\end{align}
Here, $\rho_L$ and $\rho_R$ are the so called left and right semi-radii of the Weyl disk.
The inequality \eqref{W10} may be rewritten in the form of the Weyl disk parametrization
of the values $\phi(r,\la)$ (similar, e.g., to the parametrization \cite[(2.19)]{FKRS}
for Dirac systems):
\begin{align} & \label{Z8}
\phi(r,\la)=\rho_L(r,\la)\om(r,\la) \rho_R(r,\la)-\mf{A}_{22}(r,\la)^{-1}\mf{A}_{21}(r,\la) \quad (\om^*\om\leq I_p),
\end{align}
where $\om(r,\la)$ are $p\times p$ matrices and  $\phi \in \cln_r$. Recall that the matrix inequality $B_2\geq B_1\geq 0$
yields $B_2^{1/2}\geq B_1^{1/2}$ (see, e.g., \cite{Bel}).
Hence, in view of \eqref{Z4}--\eqref{Z7}
the left and right semi-radii are non-increasing.

By $L^2(H)$ we denote the space of vector functions on $\BR_+$ with the scalar product
$$(f_1,f_2)_H=\int_0^{\infty}f_2(x)^*H(x)f_1(x)dx.$$
\begin{Pn} \label{PnH} Let $H(x)$ $(x\geq 0)$ be the Hamiltonian of  a canonical system.
Then, there is $\vp(\la)$, which  satisfies \eqref{W14} for this system. If \eqref{W14} holds, the columns of $W(x,\la)\begin{bmatrix}I_p \\ \vp(\la)\end{bmatrix}$
belong $L^2(H)$, that is,
\begin{align} & \label{W15}
\int_0^{\infty}\begin{bmatrix}I_p & \vp(\la)^*\end{bmatrix}W(x,\la)^*H(x)W(x,\la)\begin{bmatrix}I_p \\ \vp(\la)\end{bmatrix}dx<\infty \quad (\la \in \BC_+).
\end{align}
\end{Pn}
\begin{proof}. We already proved that  $\bigcap_{r>0}\cln(r)$ is non-empty. Moreover, in view of \eqref{W10}, for any $\vp$ satisfying \eqref{W14} and any $r>0$ we have
\begin{align} & \label{W16}
\begin{bmatrix}I_p & \vp(\la)^*\end{bmatrix}\mathfrak{A}(r,\la)\begin{bmatrix}I_p \\ \vp(\la)\end{bmatrix}\geq 0.
\end{align}
Taking into account \eqref{W1} and \eqref{W16}, we derive
\begin{align} &\nn
\int_0^{r}\begin{bmatrix}I_p & \vp(\la)^*\end{bmatrix}W(x,\la)^*H(x)W(x,\la)\begin{bmatrix}I_p \\ \vp(\la)\end{bmatrix}dx
\\ &  \label{W17}
\leq \frac{\I}{\la - \ov{\la}}
\begin{bmatrix}I_p & \vp(\la)^*\end{bmatrix}j\begin{bmatrix}I_p \\ \vp(\la)\end{bmatrix}\leq \frac{\I}{\la - \ov{\la}}\, I_p,
\end{align}
and \eqref{W15} follows.
\end{proof}
\begin{Dn}\label{DnW} Holomorphic $($in $\BC_+)$ $p\times p$ matrix functions $\vp(\la)$, such that the inequality \eqref{W15} holds,
are called Weyl--Titchmarsh $($Weyl$)$ functions of the canonical system  \eqref{1.3} on $[0, \, \infty)$.
\end{Dn}
Proposition \ref{PnH} implies that Weyl function always exists.
\section{Canonical systems and matrix string and Schr\"odinger equations: interconnections} \label{String}
\setcounter{equation}{0}
{\bf 1.} In view of \eqref{1.3-}, canonical systems \eqref{1.3} with Hamiltonians $H(x)$ of the form \eqref{I1} may be transformed into  systems \eqref{1.1} with Hamiltonians $\clh$:
\begin{align} & \label{B1}
\clw^{\prime}(x,\la)=\I \la J\clh(x)\clw(x,\la), \quad \clh=\vt(x)^*\vt(x), \quad \vt(x)J\vt(x)^*=0,
\end{align}
using the transformation
\begin{align} & \label{B2}
\clw=\Theta w, \quad \clh=\Theta H\Theta^*, \quad \vt(x)=\b\Theta^*.
\end{align}
Clearly, the inverse transformation works as well, that is,  systems \eqref{1.3}, \eqref{I1} and systems \eqref{B1} are equivalent.

It will be convenient to  repeat here the transformation (from \cite[Ch. 4]{SaL2-} or \cite[Section 11.1]{SaL2}) of the system \eqref{B1} into
the matrix string equation. We partition $p \times 2p$ matrix function $\vt(x)$ into $p\times p$ blocks $\vt(x)=\begin{bmatrix}\vt_1(x) & \vt_2(x)\end{bmatrix}$.
We assume that
$\det (\vt_1(x))\not=0,$ 
and we also require that $\vt_1(x)^{-1}\vt_2(x)$ is absolutely continuous and its derivative is invertible.
We set
\begin{align} & \label{B4}
\cly(x,\la)=\vt(x)\clw(x,\la), \quad \clz(x,\la)=\vt_1(x)^{-1}\vt(x)\clw(x,\la),
\\  & \label{B4+}
\vk(x):=\Big(\I\big(\vt_1(x)^{-1}\vt_2(x)\big)^{\prime}\Big)^{-1}=\vk(x)^*.
\end{align}
The self-adjointness of $\vk(x)$ above follows from the last equality in \eqref{B1}. 
According to \eqref{B1} and \eqref{B4}, we have
\begin{align} & \label{B5}
\clw^{\prime}(x,\la)=\I \la J \vt(x)^*\cly(x,\la), \\
& \label{B6}
\clz^{\prime}(x,\la)=\big(\vt_1(x)^{-1}\vt(x)\big)^{\prime}\clw(x,\la)=
\begin{bmatrix}0 &  \big(\vt_1(x)^{-1}\vt_2(x)\big)^{\prime}\end{bmatrix}\clw(x,\la).
\end{align}
Finally, taking into account \eqref{B4+}--\eqref{B6}, we see that $\clz(x,\la)$ satisfies the matrix string equation
\begin{align} & \label{B7}
\frac{d}{dx}\Big(\vk(x)\frac{d}{dx}\clz(x,\la)\Big)=\la \vt_1(x)^*\cly(x,\la)=\la \om(x)\clz(x,\la), \\
& \label{B8}
 \om (x):=\vt_1(x)^*\vt_1(x)>0.
\end{align}
{\bf 2.} Now, consider the matrix Schr\"odinger equation
\begin{align} & \label{B9}
-\clz^{\prime\prime}(x,\la)+u(x)\clz(x,\la)=\la \clz(x,\la) \quad \big(u(x)=u(x)^*\big),
\end{align}
where $u$ is a $p\times p$ matrix function.    The transformation of \eqref{B9} into the canonical system of the form \eqref{B1},
such that
\begin{align} & \label{B10}
\vt^{\prime\prime}(x)=u(x)\vt(x),
\end{align}
and  $\vt(x)$ is normalized at $x=0$ by
\begin{align} & \label{B12}
B(0)=\Theta_1:=\frac{1}{\sqrt{2}}\begin{bmatrix} \I I_p & I_p \\ \I I_p & -I_p\end{bmatrix},
\end{align}
where 
\begin{align} & \label{B11}
B(x):=\begin{bmatrix} \vt(x) \\ \vt^{\prime}(x)\end{bmatrix},
\end{align}
is described in \cite[Section 11.2]{SaL2}.
The interconnections between the spectral theories of systems \eqref{B1}, \eqref{B10} and equations \eqref{B9}
are also studied there. 
It is easily checked (see also \cite{Rem} for the case $p=1$) that the above-mentioned transformation in  \cite[Section 11.2]{SaL2}  works in the opposite direction as well.

Namely, {\it starting from  the canonical system \eqref{B1}, \eqref{B10}, \eqref{B12} one comes to the Schr\"odinger equation \eqref{B9}}. 
Indeed, according to \cite[(2.10)]{SaL2}, we have
\begin{align} & \label{B13}
B(0)^*J_1B(0)=\Theta_1^*J_1\Theta_1=J, \quad J_1:=\I \begin{bmatrix} 0 & -I_p \\ I_p & 0\end{bmatrix},
\end{align}
where $J$ is given in \eqref{1.1}. Moreover, the equalities \eqref{B10} and \eqref{B11} yield
\begin{align} & \label{B14}
B^{\prime}(x)=\begin{bmatrix} 0  & I_p \\ u(x) & 0\end{bmatrix}B(x).
\end{align}
The relations \eqref{B13} and \eqref{B14} imply that
\begin{align} & \label{B15}
B(x)^*J_1B(x)=B(0)^*J_1B(0)=J,
\end{align}
and so
\begin{align} & \label{B16}
B(x)JB(x)^*=J_1.
\end{align}
\begin{Rk}\label{nuu}
Formula \eqref{B16} shows that the equalities
\begin{align} & \label{B17}
\vt(x)J\vt(x)^*=0, \quad \vt^{\prime}(x)J\vt(x)^*=\I I_p
\end{align}
follow from \eqref{B10}--\eqref{B11}. 
\end{Rk}
Finally, setting
\begin{align} & \label{B18}
\clz(x,\la)=\vt(x)\clw(x,\la)
\end{align}
and taking into account \eqref{B1}, \eqref{B10} and \eqref{B17}, we derive
$$\clz^{\prime\prime}=u\clz-2\la \clz+\la\clz,$$
that is, $\clz(x,\la)$ satisfies matrix Schr\"odinger equation \eqref{B9}.
Formula \eqref{B18} describes the connection between the solutions of the canonical system \eqref{B1}, \eqref{B10}, \eqref{B12} and of  the corresponding
Schr\"odinger equation \eqref{B9}.

{\bf 3.} Since matrix Schr\"odinger equations may be transformed (see \cite{SaL2}) into canonical systems satisfying
\eqref{B10}--\eqref{B11} (and by virtue of Remark \ref{nuu}), they are also equivalent to a subclass of canonical systems
\eqref{1.3} with Hamiltonians of the form \eqref{F1}.
\begin{Rk}\label{NoMS} It is easy to see that in the case of our explicit formulas \eqref{E2}, \eqref{E4}, \eqref{E5}
the matrix function $\vt(x)=\b(x)\Th^*$ satisfies \eqref{B10}, where $u=-c^2 I_p$. However,
$\wt \vt(x)=\wt \b(x)\Th^*$  does not  satisfy \eqref{B10} $($excluding, possibly, some special cases$)$.
\end{Rk}
\begin{Rk}\label{RkFS} Formulas \eqref{B4}, \eqref{B4+}, \eqref{B7}, and \eqref{B8} show that  explicit
expressions for the Hamiltonians and fundamental solutions constructed in this paper generate
explicit expressions for the matrix string equations and their solutions as well.
\end{Rk}
\section{On linear similarity to squared integration} \label{Sim}
\setcounter{equation}{0}
We will consider similarity transformations of linear integral operators $K$ in $L_2^p(0,\bT)$ $(0<\bT<\infty)$:
\begin{align} & \label{AC1}
K=\I \b(x)j\int_0^x\b(t)^*\cdot dt, \quad \b(x)j\b(x)^*\equiv 0, \quad \b^{\prime}(x)j\b(x)^*\equiv \I I_p,
\end{align}
where $\b(x)$ is  a $p \times 2p$ matrix function and
\begin{align} & \label{AC1+}
 j=\begin{bmatrix} I_{p} & 0 \\  0 & -I_{p}\end{bmatrix}.
\end{align}
Recall that the operator $A$ is introduced in \eqref{AC2}.
The class of operators $K=\int_0^x K(x,t) \, \cdot \,  dt$, which are linear similar to $A$ above, was studied (for the case of the scalar kernel function
$K(x,t)$) in the essential for our considerations  paper \cite{SaL0}. Here, we study an important special subclass \eqref{AC1} of such operators
under reduced smoothness conditions on $K(x,t)$. We include the matrix case (i.e., the case $p>1$) and present a complete proof of
the similarity result.
\begin{Tm}\label{TmSim} Let operator $K$ be given by the first equality in \eqref{AC1}, and  let $\b(x)$ satisfy the second 
and third equalities in \eqref{AC1}.
Assume that $\b(x)$ is two times differentiable and the entries of $\b^{\prime\prime}(x)$ are square-integrable,
that is, 
$\b^{\prime\prime}(x)\in L_2^{p\times 2p}(0,\bT)$.  Then,  $K$ is linear similar to $A\, :$
\begin{align} & \label{AC3}
K=VAV^{-1}, \quad V=u(x)\big(I+\int_0^x \clv(x,t)\, \cdot \, dt\big),
\end{align}
where $u(x)$ is a $p\times p$ matrix function, which is unitary $($i.e., $u^*=u^{-1})$ and absolutely continuous  on $[0,\bT]$,  and 
\begin{align} & \label{AC4}
\sup\|\clv(x,t)\|<\infty \quad (0\leq t\leq x\leq \bT).
\end{align}
\end{Tm}
\begin{proof}. In the proof, we construct an operator $V$, which satisfies theorem's conditions. This $V$ is closely related
to {\it transformation operators in inverse spectral and scattering theories}.  

Step 1. Together with $K$, we consider the operators:
\begin{align} & \label{AC5}
\breve K: =\I \b^{\prime \prime}(x)j \int_0^x\b(t)^*\, \cdot \, dt, \quad (I-\Breve K)^{-1}=\int_0^x \clr(x,t) \, \cdot \,  dt.
\end{align}
The operator $\Breve K$ has a semi-separable kernel, and so (see, e.g., \cite[Section IX.2]{GGK}) the matrix function $\clr$ in \eqref{AC5} has the form
\begin{align} & \label{AC6}
\clr(x,t)=\I \b^{\prime \prime}(x)u_1(x)u_1(t)^{-1}j \b(t)^* \quad (0\leq t\leq x),
\end{align}
where the $2p\times 2p$ matrix function $u_1$ is the normalized fundamental solution of the system
\begin{align} & \label{AC7}
u_1^{\prime}(x)=\I j \b(x)^*\b^{\prime \prime}(x)u_1(x), \quad u_1(0)=I_{2p}.
\end{align}
Introduce the $p\times p$ matrix function $g(x)$ by the equalities
\begin{align} & \label{AC8}
g(0)=I_p, \quad g^{\prime}(0)=\frac{\I}{2}\b^{\prime}(0)j\b^{\prime}(0)^*, 
\\ & \label{AC9}
g^{\prime\prime}(x)=\I (I-\Breve K)\Big(\b^{\prime \prime}(x)u_1(x)\big(j\b^{\prime}(0)^*+\frac{\I}{2} j\b(0)^*\b^{\prime}(0)j\b^{\prime}(0)^*\big)\Big),
\end{align}
where the operator $(I-\Breve K)$ on the right-hand side of \eqref{AC9} is applied columnwise to the $p\times p$ matrix function above.
Further in the proof, we study the matrix function
\begin{align} & \label{AC10}
y(x,z):=(I-z^2K)^{-1}g(x).
\end{align}
(Since $K$ has a semi-separable kernel, one can write down a more explicit expression for $y$ as well.) 
Differentiating two times both parts of the equality $(I-z^2K)y(x,z)=g(x)$ and taking into account \eqref{AC1}, \eqref{AC5},  we derive
\begin{align}\nn
y^{\prime \prime}(x,z)&=g^{\prime \prime}(x)-z^2y(x,z)+\I z^2\b^{\prime \prime}(x)j\int_0^x\b(t)^*y(t,z)dt
\\  & \label{AC11}
=g^{\prime \prime}(x)-z^2(I-\Breve K)y(x,z).
\end{align}
It also easily follows from \eqref{AC1}, \eqref{AC8}, and \eqref{AC10} that
\begin{align} & \label{AC11+}
y(0,z)=I_p, \quad y^{\prime}(0,z)=\frac{\I}{2}\b^{\prime}(0)j\b^{\prime}(0)^*.
\end{align}

Applying (columnwise) $(I-\Breve K)^{-1}$ to both parts of \eqref{AC11}, using \eqref{AC5} and  \eqref{AC9}, and integrating by parts, we obtain
\begin{align} \nn
y^{\prime \prime}(x,z)=&-\int_0^x\clr(x,t)y^{\prime \prime}(t,z)dt
\\ \nn &
+\I\b^{\prime \prime}(x)u_1(x)\big(j\b^{\prime}(0)^*+\frac{\I}{2} j\b(0)^*\b^{\prime}(0)j\b^{\prime}(0)^*\big)
-z^2 y(x,z)
\\ \nn = &
-\clr(x,t)y^{\prime}(t,z)\Big|_{0}^x+\Big(\frac{\p}{\p t}\clr(x,t)y(t,z)\Big)\Big|_{0}^x
\\ \nn & 
+\I\b^{\prime \prime}(x)u_1(x)\big(j\b^{\prime}(0)^*+\frac{\I}{2} j\b(0)^*\b^{\prime}(0)j\b^{\prime}(0)^*\big)
-z^2 y(x,z)
\\ & \label{AC12}
-\int_0^x\Big(\frac{\p^2}{\p t^2}\clr(x,t)\Big)y(t,z)dt.
\end{align}
In view of the last equality in \eqref{AC1}, we have $\big(\b^{\prime}(x)j\b(x)^*\big)^{\prime}=0$, which yields
\begin{align} & \label{AC13}
u_2(x):=\I\b^{\prime\prime }(x)j\b(x)^*=-\I\b^{\prime}(x)j\b^{\prime}(x)^*.
\end{align}
Here, $u_2$ is a $p\times p$ matrix function. Relations \eqref{AC6}, \eqref{AC7}, and \eqref{AC13} imply:
\begin{align}  \label{AC14}&
\clr(x,x)=u_2(x); \quad \frac{\p}{\p t}\clr(x,t)\Big|_{t=x}=-u_2(x)^2+u_3(x), 
\\   \label{AC15}&
u_2(x)^{\prime}= u_3(x)^*-u_3(x), \quad u_3(x):=\I \b^{\prime\prime }(x)j\b^{\prime}(x)^*;
 \\ \label{AC16}&
 \frac{\p}{\p t}\clr(x,t)\Big|_{t=0}=\I \b^{\prime\prime }(x)u_1(x)\big(j\b^{\prime}(0)^*
 -j\b(0)^*u_2(0)\big);
  \\ \label{AC16+}&
\clr(x,0)=\I \b^{\prime \prime}(x)u_1(x)j \b(0)^*  .
\end{align}
Taking into account \eqref{AC11+} and \eqref{AC13}--\eqref{AC16}, we rewrite \eqref{AC12} in the form:
\begin{align} \nn
y^{\prime \prime}(x,z)+z^2 y(x,z)=&-u_2(x)y^{\prime}(x,z)-\big(u_2(x)^2-u_3(x)\big)y(x,z)
\\ \label{AC17} &
-\int_0^x\Big(\frac{\p^2}{\p t^2}\clr(x,t)\Big)y(t,z)dt,
\end{align}
where
\begin{align} & \label{AC18}
\frac{\p^2}{\p t^2}\clr(x,t)=h_1(x)h_2(t), \quad h_1(x)=\I \b^{\prime\prime}(x)u_1(x)\in L_2^{p\times 2p}(0,\bT),
\\ & \nn
h_2(t)=u_1(t)^{-1}\big(j\b(t)^*u_2(t)^2-j\b(t)^*u_3(t)^*-j\b^{\prime}(t)^*u_2(t)+j\b^{\prime\prime}(t)^*\big),
\\ & \label{AC20}
h_1\in L_2^{p\times 2p}(0,\bT), \quad h_2\in L_2^{2p\times p}(0,\bT).
\end{align}
We introduce the $p\times p$ matrix functions $y_1(x,z)$ and $u(x)$ by the equalities
\begin{align} & \label{AC21}
y_1(x,z)=u(x)^{-1}y(x,z); \quad u^{\prime}(x)=-\frac{1}{2}u_2(x)u(x), \quad u(0)=I_p.
\end{align}
Since $u_2^*=-u_2$ and $u(0)=I_p$, we obtain $u^*=u^{-1}$. Thus, using
formulas \eqref{AC15}, \eqref{AC18}, and \eqref{AC21}, we rewrite \eqref{AC17} in the form
\begin{align} & \label{AC22}
y_1^{\prime\prime}(x,z)+z^2y_1(x,z)=u_4(x)y_1(x,z)-u(x)^{*}h_1(x)\int_0^xh_2(t)u(t)y_1(t,z)dt,
\\  & \label{AC23}
u_4(x):=u(x)^{*}\Big(\frac{1}{2}\big(u_3(x)+u_3(x)^*\big)-\frac{3}{4}u_2(x)^2\Big)u(x).
\end{align}
In view of \eqref{AC11+} and \eqref{AC21}, the initial conditions for $y_1$ take the form
\begin{align} & \label{AC24}
y_1(0,z)=I_p, \quad y_1^{\prime}(0,z)=y^{\prime}(0,z)-u^{\prime}(0)=0.
\end{align}

Step 2. Let us construct the solution of system \eqref{AC22} with the initial conditions \eqref{AC24} as a series
\begin{align}  \label{AC25}
y_1(x,z)=&\sum_{k=0}^{\infty}\psi_k(x,z), \quad \psi_0(x,z)=\cos(zx)I_p,
\\ \nn
\psi_k(x,z):=&\int_0^x\cos\big(z(x-t)\big)\int_0^t\Big(u_4(s)\psi_{k-1}(s,z)
\\ & \label{AC26}
-\int_0^s\clf(s,\eta))\psi_{k-1}(\eta,z)d\eta\Big)dsdt \quad (k\geq 1),
\\ \label{AC27}
\clf(x,t):=&u(x)^{*}h_1(x)h_2(t)u(t).
\end{align}
Clearly, for $k\geq 1$ we have
\begin{align} & \label{AC28}
\psi_k^{\prime\prime}(x,z)=-z^2\psi_k(x,z)+u_4(x,z)\psi_{k-1}(x,z)-\int_0^x\clf(x,\eta)\psi_{k-1}(\eta,z)d\eta.
\end{align}
Thus, if the corresponding series converge, the matrix function $y_1$ given by \eqref{AC25}--\eqref{AC27}
satisfies \eqref{AC22}. Convergences follow from the representation
\begin{align} & \label{AC29}
\psi_k(x,z)=\int_0^x\cos (z \zeta)\clv_k(x,\zeta)d\zeta \quad (k\geq 1),
\end{align}
which is proved by induction. Indeed, setting $k=1$ in \eqref{AC26}, taking into account that $\psi_0(s,z)=\cos(zs)$
and 
\begin{align} & \label{AC30}
\cos\big(z(x-t)\big)\cos(zs)=\frac{1}{2}\Big(\cos\big(z(x-t-s)\big)+\cos\big(z(x-t+s)\big)\Big),
\\  & \label{AC30+}\cos\big(z(x-t)\big)\cos(z\eta)=\frac{1}{2}\Big(\cos\big(z(x-t-\eta)\big)+\cos\big(z(x-t+\eta)\big)\Big),
\end{align}
and changing variables  and order of integration, we derive:
\begin{align} \label{AC31}
\psi_1(x,z)=&\int_0^x\cos (z \zeta)\clv_1(x,\zeta)d\zeta,
\\ \nn
\clv_1(x,\zeta)=&\frac{1}{2}\left(\int_0^{(x+\zeta)/2}u_4(t)dt + \int_0^{(x-\zeta)/2}u_4(t)dt -\int_{(x+\zeta)/2}^x\Breve \clf(t, x-t+\zeta)dt\right.
\\  & \label{AC32}
\left. 
- \int_{(x-\zeta)/2}^{x-\zeta}\Breve \clf(t, x-t-\zeta)dt-\int_{(x-\zeta)}^x\Breve \clf(t, \zeta+t-x)dt\right),
\\   \label{AC33}
\Breve \clf(t, \eta):=&\int_{\eta}^t\clf(s,\eta)ds.
\end{align}
Assuming that \eqref{AC29} holds for $k-1$ and $\psi_k$ is given by \eqref{AC26}, 
we use a similar procedure (i.e., formulas of the  \eqref{AC30} type, change of variables
and change of order of integration) and obtain \eqref{AC29} for $k$, where
\begin{align} \nn
2\clv_k(x,\zeta)=&\int_{x-\zeta}^x\int_{\zeta+t-x}^tu_4(s)\clv_{k-1}(s,\zeta+t-x)dsdt
\\ \nn &
+\int_{(x+\zeta)/2}^x\int_{\zeta+x-t}^t u_4(s)\clv_{k-1}(s,\zeta+x-t)dsdt
\\  \label{AC34} &
+\int_{(x-\zeta)/2}^{x-\zeta}\int_{x-t-\zeta}^t u_4(s)\clv_{k-1}(s,x-t-\zeta)dsdt
\\ \nn &
-\int_{x-\zeta}^{x}\int_{\zeta+t-x}^t \int_{\zeta+t-x}^s \clf(s,\eta)\clv_{k-1}(\eta,\zeta+t-x)d\eta dsdt
\\ \nn &
-\int_{(x+\zeta)/2}^{x}\int_{\zeta+x-t}^t \int_{\zeta+x-t}^s \clf(s,\eta)\clv_{k-1}(\eta,\zeta+x-t)d\eta dsdt
\\ \nn &
-\int_{(x-\zeta)/2}^{x-\zeta}\int_{x-t-\zeta}^t \int_{x-t-\zeta}^s \clf(s,\eta)\clv_{k-1}(\eta,x-t-\zeta)d\eta dsdt.
\end{align}
Thus, the representation \eqref{AC29} is proved. Moreover, one can choose such $C(T)=C>0$ that
\begin{align} & \label{AC35}
 \int_0^\bT\|h_k(t)\|dt \leq C \quad (k=1,2),  \quad \int_0^\bT\|u_4(t)\|dt \leq C^2,
\\ & \label{AC36}
\sup_{0\leq \zeta \leq x \leq \bT} \|\clv_1(x,\zeta)\|\leq C.
\end{align}
We assume that the inequalities \eqref{AC35} and \eqref{AC36} hold. In particular, the inequality
\begin{align} & \label{AC37}
\|\clv_k(x,\zeta)\|\leq \frac{(3C^2)^{k-1}}{(k-1)!}Cx^{k-1} \quad (k\geq 1)
\end{align}
is fulfilled for $k=1$. If \eqref{AC37} is valid for $\clv_{k-1}$, relations \eqref{AC34}--\eqref{AC36} imply that \eqref{AC37} is valid for $\clv_{k}$.
Hence, \eqref{AC37} is proved.  Therefore, the series $\sum_{k=1}^{\infty}\|\clv_k(x,\zeta)\|$ is convergent.  Thus, 
the series in \eqref{AC25} converges as well, and (in view of \eqref{AC25}, \eqref{AC29}, \eqref{AC37}) we have
\begin{align} & \label{AC38}
y_1(x,z)=\cos(zx)I_p+ \int_0^x\cos (z \zeta)\clv(x,\zeta)d\zeta, \\
& \label{AC39}
 \clv(x,\ze):=\sum_{k=1}^{\infty}\clv_k(x,\ze), \quad \sup_{0\leq \zeta \leq x \leq \bT} \|\clv(x,\zeta)\|< \infty.
\end{align}
It is immediate from \eqref{AC38} and \eqref{AC39} that $y_1(0,z)=I_p$.  In order to calculate the initial value $y_1^{\prime}(0,z)$,
we differentiate both sides of \eqref{AC26} and obtain 
\begin{align} \nn
\psi_k^{\prime}(x,z)=&-z\int_0^x\sin\big(z(x-t)\big)\int_0^t\Big(u_4(s)\psi_{k-1}(s,z)
\\  \label{AC40} &
-\int_0^s\clf(s,\eta))\psi_{k-1}(\eta,z)d\eta\Big)dsdt 
\\ & \nn
+\int_0^x\Big(u_4(s)\psi_{k-1}(s,z)
-\int_0^s\clf(s,\eta))\psi_{k-1}(\eta,z)d\eta\Big)ds \quad (k\geq 1).
\end{align}
Now, the equality $y_1^{\prime}(0,z)=0$ easily follows from \eqref{AC25} and \eqref{AC40}.
Summing up, we have shown (in this step of the proof) that
$y_1$ of the form \eqref{AC38} satisfies \eqref{AC22} and \eqref{AC24}.

Step 3. Next, we show that the solution $y_1$ of  \eqref{AC22}, \eqref{AC24} is unique.
Multiplying both parts of \eqref{AC22} by the operator $A$ (given in \eqref{AC2}) and using
\eqref{AC24}, we derive
\begin{align} & \label{AC41}
By_1=I_p+z^2 Ay_1, \quad B=I-\int_0^x\clb(x,t)\cdot dt:=I+Au_4-A\int_0^x\clf(x,t)\cdot dt.
\end{align}
From \eqref{AC2}, \eqref{AC27}, and \eqref{AC41}, after simple calculations we obtain
\begin{align} & \label{AC42}
\clb(x,t)=\begin{bmatrix}I_p & x I_p & u_7(x) \end{bmatrix}\begin{bmatrix}u_5(t)-tu_4(t) \\ u_4(t)+u_6(t) \\ h_2(t)u(t) \end{bmatrix}, \\
\nn &
u_5(t):=-\int_0^tsu(s)^*h_1(s)ds \, h_2(t)u(t), \quad u_6(t):=\int_0^t u(s)^*h_1(s)ds \, h_2(t)u(t),
\\ & \label{AC43}
u_7(x):=\int_0^x (s-x)u(s)^*h_1(s)ds.
\end{align}
Here, $u_4$ is given in \eqref{AC23} and the following transformation is used:
\begin{align}  \label{AC44}
\int_0^x(t-x)\int_0^t\clf(t,s)\cdot dsdt=&\int_0^x(s-x)\int_t^x\clf(s,t)ds\cdot dt
\\ \nn
=&\int_0^x\int_0^x(s-x)u(s)^*h_1(s)ds\, h_2(t)u(t)\cdot dt
\\ \nn &
 -\int_0^x\int_0^t(s-x)u(s)^*h_1(s)ds\, h_2(t)u(t)\cdot dt.
\end{align}
Since $B$ is a triangular operator and the integral part of $B$ has a semi-separable kernel,
it easily follows (see, e.g., \cite[Section IX.2]{GGK}) that $B$ is invertible and $B^{-1}$ is a bounded operator.
(In fact, the integral part of $B$ is a Volterra operator from Hilbert-Schmidt class and $B^{-1}-I$ is again a triangular  Volterra operator with a semi-separable kernel.) 
Thus, we rewrite
\eqref{AC41} as
\begin{align} & \label{AC45+}
y_1=(I-z^2B^{-1}A)^{-1}B^{-1}I_p.\end{align}
Now, it is easy to see that $y_1$ is unique.  Recall that this unique solution admits representation
\eqref{AC3}, and so, taking into account \eqref{AC21}, we obtain
\begin{align} & \label{AC46+}
y(x,z)=V\cos(zx)I_p,
\end{align}
where $V$ is given by the second equality in \eqref{AC3}. One easily checks that
\begin{align} & \label{AC47+}
(I-z^2A)^{-1}I_p = \cos(zx) I_p.
\end{align}
In view of \eqref{AC10}, \eqref{AC46+}, and \eqref{AC47+}, we have
\begin{align} & \label{AC48+}
(I-z^2K)^{-1}g=V(I-z^2A)^{-1}I_p.
\end{align}
Presenting the resolvents in both parts of \eqref{AC48+} as series, we rewrite \eqref{AC48+}
in the form $K^ng=VA^n I_p$. In particular, setting $n=0$, we derive $g=VI_p$. The substitution $g=VI_p$
into $K^ng=VA^n I_p$ yields
\begin{align} & \label{AC49+}
K^nVI_p=VA^n I_p \quad (n\geq 0).
\end{align}
It follows that
\begin{align} & \label{AC50+}
(KV)A^nI_p=K(VA^n I_p)=K^{n+1}VI_p=VA^{n+1}I_p=(VA)A^nI_p.
\end{align}
One can easily  see (using, e.g., Weierstrass approximation theorem) that the closed linear
span of the columns of the matrix functions $A^nI_p$ ($n\geq 0$) coincides with $L_2^p(0,\bT)$.
Therefore, \eqref{AC50+} implies $KV=VA$, and \eqref{AC3} follows. The required properties of
$u$ and $\clv$ have already been proved.
\end{proof}
\begin{Rk}\label{RkLT} It is important  for the study of the canonical systems on the semi-axis $[0, \infty)$
that, according to \eqref{AC25}--\eqref{AC27}, \eqref{AC29}, and \eqref{AC39}, the matrix function
$\clv(x,\zeta)$ in the domain $0\leq \zeta\leq x\leq \ell$ is uniquely determined by
$\b(x)$ on $[0,\ell]$ $($and does not depend on the choice of $\b(x)$ for $\ell<x<\bT$ and the choice of $\bT\geq \ell)$.
\end{Rk}

{\bf Acknowledgments}  {This research    was supported by the
Austrian Science Fund (FWF) under Grant  No. P29177.}

\begin{flushright}
Alexander Sakhnovich \\
Faculty of Mathematics,
University
of
Vienna, \\
Oskar-Morgenstern-Platz 1, A-1090 Vienna,
Austria, \\
e-mail: {\tt oleksandr.sakhnovych@univie.ac.at}

\end{flushright}


\begin{thebibliography}{AGKS}

\bibitem{ArD}
D.Z. Arov and H. Dym, {\it Bitangential direct and inverse problems for systems of integral and differential equations,} Cambridge University Press, Cambridge, 2012.

\bibitem{Beli}
M.I. Belishev and V.S.  Mikhailov, 
{\it Inverse problem for a one-dimensional dynamical Dirac system (BC-method),}
Inverse Problems {\bf 30} (2014),  Art. 125013.

\bibitem{Bel}
R. Bellman,  {\it Some inequalities for the square root of a positive definite matrix,} Linear Algebra Appl. {\bf 1} (1968), 321--324.

\bibitem{BoDo}
M. Bohner and  O. Do\v{s}l\'{y}, {\it Oscillation of symplectic dynamic systems,} ANZIAM J. {\bf 46} (2004),  17--32. 

\bibitem{BreS}
J. Breuer, E. Ryckman, and B. Simon, {\it Equality of the spectral and dynamical definitions of reflection,} Comm. Math. Phys. {\bf 295} (2010),  531--550.

\bibitem{Ci}
{J.L. Cieslinski,}  
{\it Algebraic construction of the Darboux matrix revisited},
 {J. Phys. A} {\bf 42} (2009), Art. 404003.

 \bibitem{CoIv} 
A. Constantin and R. Ivanov,  {\it Dressing method for the Degasperis-Procesi equatio}, {Stud. Appl. Math.} {\bf 138} (2017),  205--226.  


\bibitem{dBr}
L.  de Branges,  {\it Hilbert spaces of entire functions,} Prentice-Hall,  Englewood Cliffs, N.J., 1968.



\bibitem{DoO}
O. Do\v{s}l\'{y}, J. Elyseeva, and R. \v{S}imon Hilscher,  {\it Symplectic difference systems: oscillation and spectral theory,} Birkh\"{a}user/Springer, Cham, 2019.

\bibitem{DoH}
O. Do\v{s}l\'{y} and R. \v{S}imon Hilscher,
{\it Disconjugacy, transformations and quadratic functionals for symplectic dynamic systems on time scales,} J. Differ. Equations Appl. {\bf 7} (2001), no. 2, 265--295. 

\bibitem{EK}
J. Eckhardt and A. Kostenko,  {\it The inverse spectral problem for indefinite strings}, Invent. Math. {\bf 204} (2016), 939--977. 

\bibitem{EGNST}
J. Eckhardt, F. Gesztesy, R. Nichols, A. Sakhnovich, and G. Teschl,  {\it Inverse spectral problems for Schr\"odinger-type operators with distributional matrix-valued 
potentials,}    Differential Integral Equations {\bf 28} (2015),  505--522.

\bibitem{EKT}
J. Eckhardt, A. Kostenko, and G. Teschl, 
{\it Spectral asymptotics for canonical systems,} J. Reine Angew. Math. {\bf 736} (2018), 285--315.

\bibitem{FKRS}
B. Fritzsche, B. Kirstein, I.Ya. Roitberg, and A.L. Sakhnovich,  
{\it Weyl theory and explicit solutions of direct and inverse problems for Dirac system with a rectangular matrix potential},  Oper. Matrices {\bf 7} (2013),  183--196. 


\bibitem{FKS}
B. Fritzsche, B. Kirstein, and A.L. Sakhnovich,  
{\it  Weyl functions of generalized Dirac systems: integral representation, the inverse problem and discrete interpolation,} J. Anal. Math. {\bf 116} (2012), 17--51.

\bibitem{Ge}
F. Gesztesy,
{\it A complete spectral characterization of the
double commutation method}, {J. Funct. Anal.},  {\bf 117} (1993), 401--446.

\bibitem{GeS}
F. Gesztesy and A. Sakhnovich,  {\it The inverse approach to Dirac-type systems based on the A-function concept,} J. Funct. Anal. 279 (2020), Art. 108609.

\bibitem{GeT}
F. Gesztesy and G. Teschl,
{\it On the double commutation method}, {Proc. Amer. Math. Soc.}, {\bf 124} (1996), 1831--1840.

\bibitem{GGK}
I. Gohberg, S. Goldberg, and M.A. Kaashoek,  {\it Classes of linear operators,} Vol. I, Birkh\"auser, Basel, 1990.


\bibitem{GKS2}
I. Gohberg, M.A. Kaashoek, and A.L. Sakhnovich, 
{\it Pseudocanonical systems with rational Weyl functions: explicit
formulas and applications}, { J. Differential Equations} {\bf 146} (1998), 375--398.


\bibitem{GoKr}
I. Gohberg  and  M.G. Krein,
\textit{Theory  and  applications  of Volterra operators  in  Hilbert  space,} Transl.  math.  monographs.  {\bf 24}, Amer. Math. Soc., Providence,  RI, 1970.

\bibitem{Gu}
C.H.~Gu, H.~Hu and  Z.~Zhou. 
\textit{Darboux transformations in integrable systems. Theory and their applications to geometry,} Springer, Dordrecht, 2005.

\bibitem{JZ}
B. Jacob, K. Morris, and H. Zwart,  {\it $C_0$-semigroups for hyperbolic partial differential equations on a one-dimensional spatial domain,} 
J. Evol. Equ. {\bf 15} (2015), 493--502.

\bibitem{JLP}
 V. Jak\v{s}i\'c, B. Landon, and A. Panati,  {\it A note on reflectionless Jacobi matrices,} Comm. Math. Phys. {\bf 332} (2014), 827--838. 

\bibitem{Kau}
D.J. Kaup, {\it Simple harmonic generation: an exact method of
solution,} Stud. Appl. Math. {\bf 59} (1978), 25--35.

\bibitem{KauS}
D.J. Kaup and H. Steudel, {\it Recent results on second harmonic
generation,} Contemporary Math. {\bf 326} (2003),  33--48.

\bibitem{KoSaTe}
{A.~Kostenko, A.~Sakhnovich, and  G.~Teschl}, 
{\it Commutation Methods for Schr\"odinger Operators with Strongly Singular Potentials,}
{Math. Nachr.} {\bf 285} (2012), 
no. 4, 392--410.

\bibitem{Krein} 
M.G. Krein,  \textit{Continuous analogues of propositions on polynomials orthogonal on the unit circle}  (Russian), 
Dokl. Akad. Nauk SSSR (N.S.) \textbf{105} (1955), 637--640. 


\bibitem{KraS}
W. Kratz and R. \v{S}imon Hilscher,
{\it Rayleigh principle for linear Hamiltonian systems without controllability,} ESAIM Control Optim. Calc. Var. {\bf 18} (2012),  501--519.

\bibitem{Langer}
H. Langer,  {\it Transfer functions and local spectral uniqueness for Sturm-Liouville operators, canonical systems and strings}, Integral Equations Operator Theory {\bf 85} 
(2016),  1--23. 



\bibitem{Mar}
V.A. Marchenko,
\textit{ Nonlinear equations and operator algebras},  D.~Reidel, Dordrecht, 1988.


\bibitem{MS}
V.B. Matveev  and M.A.  Salle,  {\it Darboux transformations and
solitons}, Springer, 
Berlin, 1991. 


\bibitem{Mog}
V. Mogilevskii, 
{\it Spectral and pseudospectral functions of Hamiltonian systems: development of the results by Arov--Dym and Sakhnovich,}
Methods Funct. Anal. Topology {\bf 21} (2015), no. 4, 370--402. 

\bibitem{Rem}
C. Remling, {\it Spectral theory of canonical systems}, 
De Gruyter, Berlin, 2018. 

\bibitem{Rom}
R. Romanov,  {\it Order problem for canonical systems and a conjecture of Valent,} Trans. Amer. Math. Soc. {\bf 369} (2017),  1061--1078. 

\bibitem{RW}
R. Romanov and H. Woracek,  {\it Canonical systems with discrete spectrum,} J. Funct. Anal. {\bf 278} (2020), Art. 108318.
 
 \bibitem{Rov}
 J. Rovnyak and L.A. Sakhnovich, {\it Pseudospectral functions for canonical differential systems. II,} 
 in Oper. Theory Adv. Appl. {\bf 218}, Birkh\"auser/Springer, Basel, 2012,
pp.  583--612, 
 
\bibitem{SaA94}
A.L. Sakhnovich,  {\it Dressing procedure for solutions of
nonlinear equations and the method of operator identities,}
{Inverse problems} { \bf 10} (1994), 699--710.

 
 
\bibitem{SaA97}
A.L.~Sakhnovich,
 {\it Iterated B\"acklund--Darboux transform for
canonical systems,} { J. Funct. Anal.} {\bf 144} (1997), 359--370.

\bibitem{SaA05}
A.L.~Sakhnovich,
 {\it Second harmonic generation: Goursat problem
on the semi-strip, Weyl functions and explicit solutions,} {Inverse
Problems} {\bf 21} (2005), 703--716.

\bibitem{ALS15}
A.L. Sakhnovich, 
{\it Dynamical and spectral Dirac systems: response function and inverse problems,}
J. Math. Phys. {\bf 56} (2015), Art. 112702.


\bibitem{ALS17}
A.L. Sakhnovich, {\it Dynamical canonical systems and their explicit solutions,} Discrete Contin. Dyn. Syst. {\bf 37} (2017),  1679--1689. 

\bibitem{ALS2019} 
A.L. Sakhnovich,  {\it New ``Verblunsky-type'' coefficients of block Toeplitz and Hankel matrices and of corresponding Dirac and canonical systems,} 
J. Approx. Theory {\bf 237} (2019), 186--209.



\bibitem{ALSgrav} 
A.L. Sakhnovich, 
{\it Einstein, $\sigma$-model and Ernst-type equations and non-isospectral GBDT version of Darboux transformation,}
arXiv:2003.13024.

\bibitem{SaSaR}
A.L. Sakhnovich, L.A. Sakhnovich   and I.Ya. Roitberg,   \textit{Inverse Problems and Nonlinear Evolution Equations. 
 Solutions, Darboux Matrices and Weyl--Titchmarsh Functions}, De Gruyter,  Berlin, 2013.
 
\bibitem{SaL0} 
L.A.  Sakhnovich, {\it The spectral analysis of Volterra operators and some inverse problems} (Russian), Dokl. Akad. Nauk SSSR (N.S.) {\bf 115} (1957), 666--669. 

\bibitem{SaL1}
L.A. Sakhnovich,  {\it On  the  factorization  of  the  transfer
matrix function,} {Sov. Math. Dokl.} { \bf 17} (1976), 203--207.

\bibitem{SaL2-}
L.A. Sakhnovich. 
{\it Factorization  problems  and  operator identities,} 
{Russian Math. Surveys}   {\bf 41} (1986), 1--64.

\bibitem{LA94}
L.A. Sakhnovich, {\it The method of operator identities and problems in analysis,}  St. Petersburg Math. J. {\bf 5} (1994),  1--69.


\bibitem{SaL2}
L.A. Sakhnovich,  {\it Spectral theory of canonical differential
systems, method of operator identities,}  Birkh\"auser, Basel, 1999.

\bibitem{SaL18}
L.A. Sakhnovich, {\it Dirac equation: the stationary and dynamical scattering problems,} in  Oper. Theory Adv. Appl. {\bf 263}, pp. 407--424,
 Birkh\"auser/Springer, Cham, 2018. 

\bibitem{Su}
M. Suzuki, {\it An inverse problem for a class of canonical systems and its applications to self-reciprocal polynomials,} J. Anal. Math. {\bf 136} (2018),  273--340.

\bibitem{Wor}
H. Woracek, {\it Asymptotics of eigenvalues for a class of singular Krein strings}, Collect. Math. {\bf 66} (2015), 469--479.

\bibitem{ZM}
V.E. Zakharov and A.V. Mikhailov, {\it On the integrability of
classical spinor models in two-dimensional space-time,} {Comm.
Math. Phys.} {\bf 74} (1980), 21--40. 


\end{thebibliography}
\end{document}